\documentclass[12pt]{article}
\usepackage{amssymb}
\usepackage{amsmath}
\usepackage{amsfonts}
\begin{document}
\newcommand{\dm}{{\mathbb M}}
\newcommand{\lfoi}{\mathcal F ^{\,(\,0,\,\infty)}}
\newcommand{\lfci}{\mathcal F ^{\,(\,c,\,\infty)}}
\newcommand{\lfii}{\mathcal F ^{\,(\,\infty,\,\infty)}}
\newcommand{\fdm}{\widehat{\mathbb M}}
\newcommand{\cfdm}{\widehat{\mathcal M}_{\infty}}
\newcommand{\wat}{\mathbb W_t}
\newcommand{\vs}{\vspace{5mm}}
\renewcommand{\labelenumi}{\roman{enumi})}
\renewcommand{\baselinestretch}{1.10}
\small\normalsize
\begin{center}
{\Large\sc Microlocalization and stationary phase.}\\ $\ $ \\
{by Ricardo Garc\'{\i}a L\'opez \footnote{Partially supported by
the DGCYT, BFM2002-01240.}}
\end{center}\vs

\noindent {\bf 0. Introduction} \vs

In this paper we make an attempt to define analogues of the local
Fourier transformations defined by G. Laumon for $\ell$-adic
sheaves in the case of formal differential systems defined over a
field $K$ of characteristic zero (cf. \cite{Lau}). Our main goal
is to prove a stationary phase formula expressing the formal germ
at infinity of the Fourier transform of a holonomic $K[t]\langle
\partial_t\rangle$-module $\dm$ in terms of the formal germs defined
by $\dm$ at its singular points. When $K$ is the field $\mathbb C$
of complex numbers and the module $\dm$ is of exponential type (in
Malgrange's sense, see \cite{Mal2}), this was done by B.
Malgrange, (see loc. cit. and \cite{Mal1}). In this case the only
transformation needed is the one given by the microlocalization
functor (which corresponds to the transformation labelled
$(0,\,\infty')$ by Laumon, this is probably known). The main point
of section 1 below is to treat the case of a $K[t]\langle
\partial_t\rangle$-module with arbitrary slopes at infinity. In
order to do this, one is forced to introduce another type of
microlocalization (corresponding to Laumon's $(\infty, \infty')$
transformation) to keep track of the contribution coming from the
germ at infinity defined by $\dm$. A trascendental construction of
this $(\infty, \infty')$ transformation was explained by B.
Malgrange to the author, the construction we give here is
algebraic, although we cannot avoid the use of some trascendental
arguments in the proof. When the base field $K$ is the field of
complex numbers, we establish a $1$-Gevrey variant of the formal
stationary phase formula and we use it to give a decomposition
theorem for meromorphic connections (the decompositions obtained
are much rougher than the decomposition according to formal
slopes, but they hold at the $s$-Gevrey level, $s>0$). In section
2 we define an analogue of Laumon's $(\infty,0')$ local Fourier
transform and we use it for the study of the singularity at zero
of the Fourier transform of $\dm$ and to establish a long exact
sequence of vanishing cycles. In section 3 we make a modest
attempt to transpose part of the above constructions into the
$p$-adic setting. We define a ring of $p$-adic microdifferential
operators of finite order, we prove a division theorem for them
and we show that, in some cases, the corresponding
microlocalization functor has the relation one would expect with
the $p$-adic Fourier transform (in a sense which is made precise
in the introduction to section 3).\vs

We will use the formalism of formal slopes and a few other results
from $\mathcal D$-module theory for which to refer e.g. to
\cite{Sa0}, \cite{Sa}. Our proof of the formal stationary phase
formula follows the leitfaden of the one given by C. Sabbah in
\cite{Sa} for modules with regular singularities. I thank C.
Sabbah for his corrections and remarks on a previous version of
these notes. Part of this work was done while the author was
invited to the University of Minnesota, I thank G. Lyubeznik and
S. Sperber for their hospitality. I thank W. Messing and G.
Christol for their useful remarks.\vs

\noindent {\bf 1. Formal stationary phase} \vs

\noindent We will use the following notations:
\begin{enumerate}
\item Unless otherwise stated, all modules over a non-commutative
ring will be left modules. We denote by $\dm$ a holonomic module
over the Weyl algebra $\wat= K[t]\langle\partial_{t}\rangle$ (that
is, we assume that for all elements $m\in \dm$ there is an
operator $P\in\wat -\{ 0\}$ such that $P\cdot m=0$). The rank of
$\dm$ is defined as $\mbox{rank}(\dm):=\dim_{K (t)}K
(t)\otimes_{K[t]}\dm$. Let $\overline{K}$ be an algebraic closure
of $K$. There is a maximal Zariski open subset $U\subset \mathbb
A^1_{\overline{K}}$ such that the restriction of $\dm$ to $U$ is
of finite type over the ring of regular functions on $U$, by
definition the set of singular points of $\dm$ is
$Sing(\dm)=\mathbb A^1_{\overline{K}}-U$. We will assume the
points of $Sing(\dm)$ are in $K$.
\item  The Fourier antiinvolution is the morphism of $K$-algebras
$\wat\to\mathbb W_{\eta}$ given by $t\mapsto-\partial_{\eta}\,,\
\partial_t\mapsto\eta$. The Fourier transform of $\dm$ is defined
as the $\mathbb W_{\eta}$-module $\fdm:=\mathbb
W_{\eta}\otimes_{\wat}\dm$, where $\mathbb W_{\eta}$ is regarded
as a right $\wat$-module via the Fourier morphism. If $m\in \dm$,
we put $\widehat{m}=1\otimes m \in \fdm$.
\item
We will set $\mathcal{K}_{\eta^{-1}} :=K[[\eta^{-1}]][\eta]$ and
we consider on this field the derivation
$\partial_{\eta^{-1}}=-\eta^2\partial_{\eta}$. If $\mathcal V$ is
a $\mathcal{K}_{\eta^{-1}}$-vector space, a connection on
$\mathcal V$ is a $K$-linear map $\nabla:\mathcal V \to \mathcal V
$ satisfying the Leibniz rule
\[
\nabla(\alpha\cdot v)=\partial_{\eta^{-1}}(\alpha)\cdot v+ \alpha
\cdot \nabla (v)\ \ \mbox{for all}\ \
\alpha\in\mathcal{K}_{\eta^{-1}}\,,\ v\in\mathcal V.
\]
\item We set $\dm_{\infty}= K[[t^{-1}]]\langle
\partial_t \rangle\otimes _{K[t^{-1}]\langle
\partial_t \rangle}\dm\,[t^{-1}]$, and for $c\in K$ we
put $t_c = t-c$, $\dm _{\,c} = K[[t_c]]\langle
\partial_{t_c} \rangle \otimes_{\wat}\dm$.
\end{enumerate}

We will consider the following rings:
\begin{enumerate}
\item {\it The ring $\lfci$ of formal microdifferential operators:}

Let $c\in K$. For $r\in \mathbb Z$, we denote by $\lfci [r]$ the
set of formal sums
\[
\sum_{i\le \,r} a_i(t_c)\,\eta^i \ \ \mbox{where}\ \ a_i(t_c) \in
K[[t_c]]\,,r\in \mathbb Z.
\]
We put $\lfci=\cup_r\,\lfci[r]$. For $P,Q \in\lfci$, their product
is defined by the formula
\[
P \cdot Q = \sum_{\alpha \geqslant\, 0}\ \frac{1}{\alpha\, !}\ \
\partial_{\eta}^{\,\alpha}P\cdot
\partial^{\,\alpha}_{t_c} Q \,
\, \in \lfci.
\]
(where the product on the right hand side is the usual,
commutative product). With this multiplication, $\lfci$ becomes a
filtered ring and $\lfci [0]$ is a subring. One has morphism of
$K$-algebras given by
\begin{eqnarray*}
\varphi^{(\,c,\,\infty\,)}: \wat &\longrightarrow& \lfci \,,\\
t\ \ &\mapsto&\ t_c+c\\
\partial_t\ \ &\mapsto&\ \eta
\end{eqnarray*}
which endows $\lfci$ with a structure of $(\wat,\wat)$-bimodule.

\item {\it The ring $\lfii$:}

For $r\in \mathbb Z$, we denote by $\lfii [r]$ the set of formal
sums
\[
\sum_{i\le \, r} \ a_i(t^{-1})\, \eta^i \ \mbox{ where } \ \
a_i(t^{-1})\in K[[t^{-1}]]\,,r\in \mathbb Z.
\]
We put $\lfii=\cup_r\,\lfii[r]$. If $P,Q\in \lfii$, their product
is given by
\[
P \ast Q = \sum_{\alpha \ge\, 0}\ \frac{1}{\alpha\, !}\ \
\partial^{\,\alpha}_{\eta}P \cdot
\partial^{\,\alpha}_{t}Q
\]
Again, $\lfii$ is a filtered ring and $\lfii[0]$ is a subring. One
has a morphism of $K$-algebras
\begin{eqnarray*}
\varphi^{(\,\infty,\,\infty\,)}: K[t^{-1}]
\langle\partial_{t}\rangle &\longrightarrow& \lfii \\
t^{-1}\ \ &\mapsto&\ t^{-1}\\
\partial_{t}\ \ &\mapsto&\ \eta
\end{eqnarray*}
(notice that on the ring $K[t^{-1}] \langle\partial_{t}\rangle$,
one has the relation $[\partial_{t}, t^{-1}]=-t^{-2}$). The
morphism $\varphi^{(\,\infty,\,\infty\,)}$ endows $\lfii$ with a
structure of
$(\,K[t^{-1}]\langle\partial_{t}\rangle\,,\,K[t^{-1}]\langle
\partial_{t}\rangle\,)$
- bimodule.
\end{enumerate}

If $P=\sum_{i\in \mathbb Z} a_i(t_c)\,\eta^i\in \lfci$, the order
of $P$ is the largest integer $r$ such that $a_r(t_c)\neq 0$, and
we define the principal symbol of $P$ as
$\sigma(P)=a_{r}(t_c)\,\eta^{r}$. We define similarly the order
and the principal symbol of an operator $P\in\lfii$. Principal
symbols are multiplicative, in the sense that $\sigma(P\cdot
Q)=\sigma(P)\cdot\sigma(Q)$.\vs

We recall next some results which are well-known for the rings
$\lfci$. The proofs for $\lfii$ follow a similar pattern (using
the fact that the graded ring associated to the filtration on
$\lfii$ is isomorphic to $K[[t^{-1},x]][x^{-1}]$), and therefore
they are omitted. \vs

\noindent (1.1)\,{\bf Division theorem} (cf. e.g. \cite[Ch.4,
2.6.]{Bjo}): {\it Let $F\in\lfci$ and assume that $\sigma(F)=t_c^m
\,b(t_c)$ where $b(0)\neq 0$. Then, for all $G\in\lfci$ there
exist unique $Q\in \lfci$ and $R_0,\dots ,R_{m-1}\in \mathcal
K_{\eta^{-1}}$ such that
\[
G= Q\cdot F + R_{m-1}\, t_c^{m-1} + \dots + R_0.
\]
The same statement holds for the ring $\lfii$, replacing $t_c$ by
$t^{-1}$.} \vs

\noindent{\bf Proposition} (cf. e.g. \cite[Ch.4, 2.1 and
2.9]{Bjo}): {\it The rings $\lfci$ and $\lfii$ are left and right
Noetherian.}\vs

\noindent{\bf Proposition:} (cf. e.g. \cite[Ch.5, \S 5]{Bjo}):
{\it $\lfci$ is a flat left and right $\wat$-module. $\lfii$ is a
flat left and right $K[t^{-1}]\langle
\partial_t\rangle$-module.}\vs

We will consider the following modules:\vs

\noindent i)\,{\it The (ordinary) microlocalization of $\dm$ at
$c\in K$}: It is defined as
\[
\lfci(\dm):= \lfci \otimes_{\wat} \dm.
\]
where $\lfci$ is viewed as a right $\wat$-module via
$\varphi^{(c,\,\infty)}$. It has a structure of $\mathcal
K_{\eta^{-1}}$-vector space with a connection given by left
multiplication by \\$\eta^2 \cdot (t_c+c)=\eta^2\cdot t$. Notice
that
\[
\lfci(\dm)\cong \lfci \otimes_{K[[t_c]]\langle
\partial_{t_c}\rangle}(K[[t_c]]\langle
\partial_{t_c}\rangle\otimes_{\wat} \dm)=\lfci
\otimes_{K[[t_c]]\langle
\partial_{t_c}\rangle} \dm _c\,,
\]
thus $\lfci(\dm)$ depends only on the formal germ $\dm_c$.\vs

\noindent ii)\,{\it The $(\,\infty\,, \infty\,)$-microlocalization
of $\dm$}\,: It is defined as
\[
\lfii(\dm):= \lfii \otimes_{K[t^{-1}]\langle \partial_t\rangle}
\dm \,[t^{-1}]\,,
\]
where $\lfii$ is viewed as a $K[t^{-1}]\langle
\partial_t\rangle$-module via the morphism
$\varphi^{(\,\infty,\,\infty\,)}$. Again, it has a structure of
$\mathcal K_{\eta^{-1}}$-vector space with a connection, defined
by
\[
\nabla(\alpha\otimes m):=\partial_{\eta^{-1}}(\alpha) \otimes m
\,+\, \eta^{2}\alpha\otimes t\cdot m\,,
\]
and it depends only on $\dm_{\infty}$, since one has
\[
\lfii(\dm):=\lfii
\otimes_{K[[t^{-1}]]\langle\partial_t\rangle}\dm_{\infty}\,.
\]

By flatness of $\varphi^{(\,0,\,\infty\,)}$ and
$\varphi^{(\,\infty,\,\infty\,)}$, both microlocalizations define
exact functors.\vs

For the proof of the formal stationary phase formula we will need
the following result:\vs

\noindent {\bf Proposition 1}:\ {\it Let\
$Q(t^{-1},\partial_t)=\sum_{v=1}^n b_v(t^{-1})\,\partial^v_t \in
K[t^{-1}]\langle
\partial_t\rangle$ be such that there is at least an index $v\in
\{1,\dots,n\}$ with $b_v(0)\neq 0$. Then there is an isomorphism
of $K[t^{-1}]\langle
\partial_t\rangle$-modules
\[
\lfii\otimes_{K[t^{-1}]\langle
\partial_t\rangle}\frac{K[t^{-1}]\langle
\partial_t\rangle}{K[t^{-1}]\langle
\partial_t\rangle\cdot Q}\cong \lfii\otimes_{K[t^{-1}]\langle
\partial_t\rangle}\frac{K[t, t^{-1}]\langle
\partial_t\rangle}{K[t, t^{-1}]\langle
\partial_t\rangle\cdot Q}\,.
\]}

{\it Proof:}\ We show first that the natural map
\[
\gimel: \frac{K[t^{-1}]\langle
\partial_t\rangle}{K[t^{-1}]\langle
\partial_t\rangle\cdot Q}\longrightarrow\frac{K[t, t^{-1}]\langle
\partial_t\rangle}{K[t, t^{-1}]\langle
\partial_t\rangle\cdot Q}
\]
is injective. Assume we have $A(t,t^{-1},\partial_t)\in
K[t,t^{-1}]\langle
\partial_t\rangle$ such that $A(t,t^{-1},\partial_t)\cdot
Q(t^{-1},\partial_t)\in K[t^{-1}]\langle \partial_t\rangle$, we
have to show that in fact $A\in K[t^{-1}]\langle
\partial_t\rangle$. Write $A=\sum_u a_u(t,t^{-1})\,\partial_t^u$
with $a_u(t,t^{-1})\in K[t,t^{-1}]$. Let $v_0$ be the largest
index with $b_{v_0}(0)\neq 0$ and let $k_0\in\mathbb N$ be the
largest exponent of $t$ appearing in the Laurent polynomials
$\{a_u\}_u$ (if $a_u\in K[t^{-1}]$ for all $u$, we are done). Let
$u_0$ be the largest index such that $a_{u_0}$ contains a monomial
$\beta t^{k_0}\neq 0$, $\beta\in K$. Set $j_0=u_0+v_0$. The
coefficient of $\partial_t^{j_0}$ in $A\cdot Q$ is
\[
\sum_{{u,v,\alpha}\atop{j_0=u+v-\alpha}}\frac{1}{\alpha\,!}\,u(u-1)\dots(u-\alpha+1)\,a_u\,
\frac{d^{\alpha}b_v}{d\,t^{\alpha}}
\]
The monomial $\beta b_{v_0}(0)t^{k_0}$ appearing in the summand
corresponding to $u=u_0, v=v_0, \alpha=0$ cannot be cancelled,
because in the other summands either $u>u_0$, and then in $a_u\,
\frac{d^{\alpha}b_v}{d\,t^{\alpha}}$ all powers of $t$ appear with
exponent strictly smaller than $k_0$\,, or else $u\leqslant u_0$,
and then $\frac{d^{\alpha}b_v}{d\,t^{\alpha}}\in t^{-1}K[t^{-1}]$,
thus the exponents of $t$ in these summands are also strictly
smaller than $k_0$. Since $A\cdot Q\in K[t^{-1}]\langle
\partial_t\rangle$, we conclude that $k_0=0$, which proves the injectivity of
$\gimel$.\vs

By flatness of $\lfii$ over $K[t^{-1}]\langle
\partial_t\rangle$, the map $\mbox{Id}_{\lfii}\otimes \gimel$ is
injective as well. In order to show that it is a surjection, we
prove first the following statement:
\begin{eqnarray*}
&\mbox{{\it Claim}: For all }\ P\in K[t,t^{-1}]\langle
\partial_t\rangle \
,\mbox{ there exists a polynomial } &\\ &p(x)\in K[x]-\{ 0\}
\mbox{ such that} \ \ \ p(\partial_t)\cdot P \in
K[t,t^{-1}]\langle
\partial_t\rangle\cdot Q + K[t^{-1}]\langle \partial_t\rangle &
\end{eqnarray*}
\noindent {\it proof of the claim:} Let us denote by
$\Omega\subseteq K[t,t^{-1}]\langle
\partial_t\rangle$ the set of differential operators
$P\in K[t,t^{-1}]\langle
\partial_t\rangle$ satisfying the condition
of the claim, notice that $\Omega$ is closed under adition. It
suffices to show that $t^i\in\Omega$ for all $i\geqslant 1$,
because if there is a $p(x)$ with $p(\partial_t)\cdot
t^i=\alpha(t,t^{-1},\partial_t)\cdot Q + \beta(t^{-1},\partial_t)$
then, multiplying on the right by $\partial_t^j$, we obtain that
$\partial_t^j\cdot t^i\in\Omega$ for all $i, j\geqslant 0$, and
then we are done. We prove $t^i\in \Omega$ by induction on
$i\geqslant 1$: By our hypothesis on $Q$ there exists a non-zero
polynomial $p(x)\in K[x]$ such that $Q=
p(\partial_t)+\beta(t^{-1},\partial_t)$ with
$\beta(t^{-1},\partial_t)\in t^{-1}K[t^{-1}]\langle
\partial_t\rangle$ (we will denote $p\,'=\frac{dp}{dx}$).
Now, multiplying this equality on the left by $t$ and using
$t\,p(\partial_t)=p(\partial_t)\,t-p\,'(\partial_t)$, case $i=1$
follows. For the induction step, assume we have
\[
p(\partial_t)\cdot t^i =\alpha(t,t^{-1},\partial_t)\cdot Q +
\beta(t^{-1},\partial_t)
\]
We have also
$p(\partial_t)\,t^{i+1}=t\,p(\partial_t)t^i+p\,'(\partial_t)\,t$,
substituting,
\[
p(\partial_t)\,t^{i+1}=t(\alpha\,Q+\beta)+p\,'(\partial_t)\,t
\]
Now from the case $i=1$ follows that we have
$t\,\beta(t^{-1},\partial_t)\in\Omega$ and
$p\,'(\partial_t)\,t\in\Omega$, so the claim is proved. Given
\[
F(t^{-1},\eta)\otimes
P(t,t^{-1},\partial_t)\in\lfii\otimes\frac{K[t, t^{-1}]\langle
\partial_t\rangle}{K[t, t^{-1}]\langle
\partial_t\rangle\cdot Q}
\]
choose $p(x)\neq 0$ such that $p(\partial_t)\cdot P\in
K[t,t^{-1}]\langle
\partial_t\rangle\cdot Q + K[t^{-1}]\langle \partial_t\rangle$.
By the division theorem $p(\eta)$ is invertible in $\lfii$, thus
we have $F\otimes P = F\cdot p(\eta)^{-1}\otimes
p(\partial_t)\,P\in \mbox{Im}\, [\,\mbox{Id}_{\lfii}\otimes
\gimel]$, and then the proposition is proved. $\Box$ \vs

\noindent {\it Definition}:\ If $\mathbb N$ is a $\mathbb W_{\eta}
$-module, its formal germ at infinity is the $\mathcal
K_{\eta^{-1}}$-vector space $\mathcal N_{\infty}=\mathcal K
_{\eta^{-1}}\otimes_{K[\eta]} \mathbb N$, which is endowed with
the connection defined by
\[
\nabla (\alpha \otimes n)=\partial_{\eta^{-1}}(\alpha)\otimes n -
\alpha \otimes \eta^2\partial_{\eta}n.
\]
The main result of this section is:\vs

\noindent {\bf Theorem}\ (formal stationary phase): {\it Let $K$
be a field of characteristic zero, let $\dm$ be a holonomic
$K[t]\langle
\partial_t \rangle$-module. Then, after a finite extension
of the base field $K$, the map
\[
\Upsilon : \cfdm \longrightarrow\bigoplus_{c\in\,
Sing\dm\cup\{\infty\}} \lfci (\dm)
\]
given by $\Upsilon (\alpha\otimes \widehat{m}) = \oplus_c \
\alpha\otimes m$ is an isomorphism of $\mathcal
K_{\eta^{-1}}$-vector spaces with connection. } \vs

\noindent {\it Proof:} The connections on the right hand side have
been chosen so that the map is a morphism of $\mathcal
K_{\eta^{-1}}$-vector spaces with connection, we have to show it
is an isomorphism. We will assume all $\wat$-modules appearing in
what follows have $K$-rational singularities (which can be
achieved after a finite extension of $K$). Consider first the
Dirac modules $\delta_c= \wat /\wat (t-c)$. It is easy to check
that we have
\begin{eqnarray*}
&&\lfii(\delta_c)=0\,,\ \mathcal F ^{\,(\,d,\,\infty)}(\delta_c)=0
\mbox{ if }d\neq c\,, \\
&&\mathcal K _{\eta^{-1}}\otimes \widehat{\delta_c}=\mathcal K
_{\eta^{-1}} = \lfci(\delta_c)\,,
\end{eqnarray*}
so the theorem follows in this case.\vs

For an arbitrary holonomic module $\mathbb M$, there is a
differential operator $P\in\wat$ and a $\wat$-module with punctual
support $\mathbb K$ so that one has an exact sequence
\[
0 \longrightarrow \mathbb K \longrightarrow \wat / \wat \cdot P
\longrightarrow \dm \longrightarrow 0\,.
\]
A holonomic module with punctual support is a finite direct sum of
Dirac's $\delta$-modules, since  both the global and the local
Fourier transforms are exact functors, it will suffice to consider
the of a quotient of $\wat$ by an operator. Moreover, given
$\mathbb M=\wat /\wat\cdot P_1$, there is an operator $P_2\in\wat$
such that one has
\[
0 \longrightarrow \mathbb K_1 \longrightarrow \wat / \wat \cdot
P_1=\mathbb M \longrightarrow \mathbb M\,[t^{-1}]=\wat/\wat \cdot
P_2\longrightarrow \mathbb K_2 \longrightarrow 0\,.
\]
where $\mathbb K_i$ are supported at zero for $i=1,2$
(\cite[4.2]{Sa0}). Thus, we will assume in what follows that $\dm
= \dm\,[t^{-1}] = \wat /\wat P $, we write
$P(t,\,\partial_t)=\sum_{i=0}^d\,
 a_i(t)\partial_t^i$ with $a_i(t)\in K[t]$. \vs

\noindent {\bf Step 1:} {\it The map $\Upsilon_c: \cfdm\to
\lfci(\dm)$, given by the composition of $\Upsilon$ with the
projection onto $\lfci(\dm)$, is exhaustive for all
$c\in\mbox{Sing}(\dm)\cup\{\infty\}$:}\vspace{4mm}

Assume first that $c=0$. Then $P(t,\eta)\in\lfoi$ has a principal
symbol of the form $\sigma(P)=t^{m_0}\,b(t)$ with $b(0)\neq 0$,
$m_0\geqslant 0$, and then we have
\[
\lfoi(\dm)=\lfoi /\lfoi\cdot P.
\]
Given $G\in\lfoi$, by the division theorem
\[
G \equiv \sum_{i=0}^{m_0-1} R_i\cdot t^i\ \ \ \ (\mbox{mod}\ \lfoi
\,P)
\]
where $R_i\in \mathcal K_{\eta^{-1}}\ (0\le i\le m_0-1)$. Then, a
preimage of the class of $G$ in $\lfoi(\dm)$ under $\Upsilon_c$ is
given by
\[
\sum_{i=0}^{m_0-1} R_i\otimes \widehat{t^i}\,\in\, \mathcal
K_{\eta^{-1}}\otimes \fdm\, ,
\]
where $\widehat{t^i}$ denotes the element $1\otimes t^i$ of $\fdm
= \mathbb W_{\eta}\otimes \dm$. The proof for arbitrary $c\in K$
is done in the same way, using the division theorem in $\lfci$
(notice that the hypothesis $\mathbb M = \mathbb M[t^{-1}]$ has
not been used up to now).\vs

We consider now the case $c=\infty$. Write
$\widehat{d}=\max_{j=1}^d\{\deg a_j(t)\}$, set
$Q(t^{-1},\partial_t)= t^{-\widehat{d}}P(t, \partial_t)$. Since
$Q$ satisfies the hypothesis of proposition 1 above, we have
\[
\frac{K[t^{-1}]\langle \partial_t\rangle}{K[t^{-1}]\langle
\partial_t\rangle \cdot Q} \cong \frac{K[t, t^{-1}]\langle \partial_t\rangle}
{K[t,t^{-1}]\langle
\partial_t\rangle \cdot Q}=\mathbb M
\]
as $K[t^{-1}]\langle \partial_t\rangle$-modules. Then we have
\[
\lfii (\dm) \cong \frac{\lfii}{\lfii \cdot Q}\,,
\]
and notice that the principal symbol of $Q$ is of the form $\sigma
(Q) = t^{\,\deg (a_d(t))-\widehat{d}}\cdot b(t^{-1})$ with
$b(0)\neq 0$. Given $G\in\lfii$, by the division theorem
\[
G\equiv \sum_{i} R_i\cdot t^{-i}\ \ \ \ (\mbox{mod}\ \lfii \,Q)
\]
where $i\in\{0,\dots,\widehat{d}-\deg (a_d(t))\}$. Since $\dm =\dm
\,[t^{-1}]$, we have that $\sum_i R_i\,\otimes\,
\widehat{t^{-i}}\in\cfdm$ is a preimage of $G$ under
$\Upsilon_{\infty}$.\vs

\noindent {\bf Step 2:} $\dim_{\mathcal K_{\eta^{-1}}}\,\cfdm=
\dim_{\mathcal K_{\eta^{-1}}}\,\left[\left(\bigoplus_{c\in Sing\,
\dm}\lfci(\dm)\right)\oplus\ \lfii(\dm_{\,\infty})\right]$.\\
\noindent Put
\[
a_d(t)=\prod_{c\in Sing\, \dm}(t-c)^{m_c}, \ \ \deg(a_d(t))=\sum\,
m_c.
\]
It follows from the existence and uniqueness of division for
$\lfci$ and $\lfii$ that we have $\dim_{\mathcal K
_{\eta^{-1}}}\lfci(\dm)=m_c$,\ $ \dim_{\mathcal K
_{\eta^{-1}}}\lfii(\dm)=\widehat{d}-\deg (a_d(t)). $ Since
$\widehat{d}=\mbox{rank}(\fdm)=\dim_{\mathcal
K_{\eta^{-1}}}\,\cfdm$, the claimed equality follows.\vs

\noindent {\bf Step 3:}{\it
\begin{itemize}
\item[a)] All slopes of $\lfoi(\dm)$ are strictly smaller than $+1$.
\item[b)] For $c\in K -\{ 0\}$, all slopes of $\lfci (\dm)$
are equal to $+1$.
\item[c)] All slopes of $\lfii(\dm)$ are strictly greater than $+1$.
\end{itemize}}

We assume first that the base field $K$ is the field $\mathbb C$
of complex numbers. Then, there is a holonomic $\mathbb
C[t]\langle \partial_t\rangle$-module $\mathbb N$ which is only
singular at $0$ and $\infty$, the singularity at infinity is
regular and for the singularity at zero we have $\mathbb C
[[t]]\langle
\partial_t\rangle \otimes_{\wat}\mathbb N \cong \mathbb C[[t]]\langle
\partial_t\rangle \otimes_{\wat}\dm$ (this is the
transcendental step in the proof, see \cite{Mal1}, it follows that
$\lfoi(\dm)\cong\lfoi(\mathbb N)$). Let $ \mathbb
L\twoheadrightarrow\mathbb N $ be a surjection where $\mathbb L$
is the quotient of $\mathbb C[t]\langle
\partial_t\rangle$ by the left ideal generated by a single
differential operator. Then we have\vs

\begin{picture}(165,80)(-80,0)

\put(0,0){$\lfoi(\mathbb L)$}

\put(0,70){$\mathcal K_{\eta^{-1}}\otimes_{\mathbb
C[\eta^{-1}]}\widehat{\mathbb L}$}

\put(140,0){$\lfoi(\mathbb N)$}

\put(140,70){$\mathcal K_{\eta^{-1}}\otimes_{\mathbb
C[\eta^{-1}]}\widehat{\mathbb N}$}

\put(70,5){$\vector(1,0){62}$}

\put(70,5){$\vector(1,0){58}$}

\put(80,75){$\vector(1,0){55}$}

\put(80,75){$\vector(1,0){51}$}

\put(30,60){$\vector(0,-1){45}$}

\put(160,60){$\vector(0,-1){45}$}

\put(38,35){$\Upsilon_0^{\mathbb L}$}

\put(165,35){$\Upsilon_0^{\mathbb N}$}

\end{picture}\vs

\noindent where we denote $\Upsilon_0^{\mathbb L}$ (respectively,
$\Upsilon_0^{\mathbb N}$) the map defined in step 1 for the module
$\mathbb L$ (resp., for $\mathbb N$) and the horizontal arrows are
surjective. By step 1, the map $\Upsilon_0^{\mathbb L}$ is onto,
and then so is $\Upsilon_0^{\mathbb N}$. But the behavior of
formal slopes under Fourier transform is well known, in particular
the slopes of $\mathcal K _{\eta^{-1}}\otimes_{\mathbb
C[\eta]}\widehat{\mathbb N}$ are strictly smaller than $+1$, and
a) follows.\vs

For b), take now a $\mathbb C[t]\langle\partial_t\rangle$-module
$\mathbb N^c$ which is only singular at $c\in\mathbb C-\{0\}$ and
at infinity, such that one has an isomorphism of formal germs
$\mathbb C [[t_c]]\langle
\partial_{t_c}\rangle\otimes\mathbb N^c\cong \mathbb C
[[t_c]]\langle \partial_{t_c}\rangle\otimes\mathbb M$,\ and the
singularity at infinity of $\mathbb N^c$ is regular. Then (as in
a) above) the corresponding map $\Upsilon^{\mathbb N^c}_c$ is
onto, and since $\mathcal K _{\eta^{-1}}\otimes_{\mathbb
C[\eta]}\widehat{\mathbb N^c}$ has only slope $+1$ at infinity, b)
follows. \vs

For the proof of c) consider first a $\mathbb
C[t]\langle\partial_t\rangle$-module $\mathbb N^0$ with a regular
singularity at zero, no other singularity and $\mathbb
N^{0}_{\infty}\cong \dm _{\infty}$ (then $\lfii(\dm)\cong
\lfii(\mathbb N^0)$). Let $\mathbb L^0\twoheadrightarrow\mathbb
N^0$ be a surjection where $\mathbb L^0$ is the quotient of
$\mathbb C[t]\langle \partial_t\rangle$ by a single differential
operator. Then $\mathbb L^0[t^{-1}]$ is given by a single
differential operator as well and, as in case a) above, the map
\[
\Upsilon^{\mathbb N^0[t^{-1}]}_{\infty}:\,\mathcal K
_{\eta^{-1}}\otimes_{\mathbb C[\eta]}\widehat{\mathbb
N^0[t^{-1}]}\longrightarrow\lfii(\mathbb
N^0[t^{-1}])\cong\lfii(\dm)
\]
is onto (the last isomorphism holds because $\dm
_{\infty}\cong\mathbb N^0_{\infty}=\mathbb N^0[t^{-1}]_{\infty}$).

Since the slopes of $K _{\eta^{-1}}\otimes_{\mathbb
C[\eta]}\widehat{\mathbb N^0[t^{-1}]}$ are either zero or strictly
greater than $+1$, the same holds for $\lfii(\dm)$.\vs

Let now $\mathbb N^1$ denote the pull back of $\mathbb N^0$ by the
translation $t\mapsto t+1$. In the same way we get a surjection
\[
\Upsilon^{\mathbb N^1[t^{-1}]}_{\infty}:\,\mathcal K
_{\eta^{-1}}\otimes_{\mathbb C[\eta]}\widehat{\mathbb
N^1[t^{-1}]}\longrightarrow\lfii(\mathbb
N^1[t^{-1}])\cong\lfii(\dm)
\]
and the slopes of $K _{\eta^{-1}}\otimes_{\mathbb
C[\eta]}\widehat{\mathbb N^1}$ are greater or equal than $+1$.
Thus all slopes of $\lfii(\dm)$ must be strictly greater than
$+1$, and then we are done when the base field is $\mathbb C$.\vs

In the general case, there is a subfield $K_1\subset K$ of finite
transcendence degree over $\mathbb Q$ such that the module $\dm$
is defined over $K_1$ and all its singular points are
$K_1$-rational. Choosing an embedding of fields
$K_1\hookrightarrow\mathbb C$, the statement follows from the
complex case treated above. \vs

If $\mathcal N$ is a $\mathcal K_{\eta^{-1}}$-vector space with
connection and $q\in\mathbb Q$, we denote by $\mathcal N_q$ its
subspace of formal slope $q$ (see e.g. \cite[5.3.1]{Sa0}), and we
denote $\mathcal N_{<q}$ its subspace of slopes strictly smaller
than $q$. If $\varphi:\mathcal L \to \mathcal N$ is a morphism, we
denote by $\varphi_{<q}: \mathcal L _{<q}\to \mathcal N_{<q}$ the
induced morphism (and similarly, we let $\varphi_{>q}$ denote the
restriction of $\varphi$ to the subspaces of slopes strictly
bigger than $q$). For $c\in\mathbb C$, let $\mathcal E^c$ denote
the one dimensional $\mathcal K_{\eta^{-1}}$- vector space with
connection given by $ \nabla (1)=c\cdot\eta^2$. It is easy to see
that, if $c\neq0$, then $\mathcal E^c$ has slope $+1$. Let
$\tau_c:K[t]\to K[t]$ the translation given by $t\mapsto t+c$. One
has an isomorphism of $\mathcal K_{\eta^{-1}}$- vector spaces with
connection \ $\mathcal E^c \otimes \lfci (\dm) \cong \lfoi
(\tau_c^{\ast}\dm)$ (from this isomorphism one can get an
alternative proof of b) above). Notice also that $\mathcal E^c
\otimes\mathcal E^{-c}\cong (\mathcal
K_{\eta^{-1}},\partial_{\eta^{-1}})$.\vs

\noindent {\bf Step 4:\ }{\it $\Upsilon$ is an isomorphism.} We
remark first that if $\mathcal M$ is $\mathcal K_{\eta^{-1}}$-
vector space with connection such that all its slopes are strictly
smaller than $+1$ then for $c\neq 0$ the twisted vector space with
connection $\mathcal E^c\otimes_{\mathcal K_{\eta^{-1}}}\mathcal
M$ has only slope $+1$ (this can be easily seen using the
structure theorem of formal meromorphic connections, \cite[Theorem
5.4.7]{Sa0}). \vs

Set $\cfdm^c:= \mathcal E^{-c} \otimes (\mathcal E^c \otimes
\cfdm)_{<1}\subset \cfdm$. If $c,d\in Sing (\dm)$ are distinct,
then the map
\[
\mbox{ id } \otimes {\Upsilon_c}_{\,\mid \cfdm^d}: \mathcal E^c
\otimes \cfdm^d \longrightarrow \mathcal E^c \otimes \lfci (\dm)
\cong \lfoi (\tau_c^{\ast}\dm)
\]
is the zero map, because its source is purely of slope $+1$ and
the target has slopes strictly smaller than $+1$. Tensoring with
$\mathcal E^{-c}$, it follows that ${\Upsilon_c}_{\,\mid \cfdm^d}:
\cfdm^d\to\lfci (\dm _{\,c})$ is the zero map as well.

Since by step one $\Upsilon_c$ is onto, also is the map $\mbox{ id
}\otimes {\Upsilon_c}:\mathcal E^c \otimes \cfdm\to\mathcal E^c
\otimes \lfci (\dm)$. Since $\mathcal E^c \otimes \lfci (\dm)\cong
\lfoi (\tau_c^{\ast}\dm)$ has slopes strictly smaller than $+1$,
the restriction
\[
(\mbox{ id }\otimes {\Upsilon_c})_{<1}: (\mathcal E^c \otimes
\cfdm)_{<1}\longrightarrow \mathcal E^c \otimes \lfci (\dm)
\]
is onto as well, thus tensoring with $\mathcal E^{-c}$ follows
that ${\Upsilon_c}_{\,\mid \cfdm^c}$ is onto.

Let $\Upsilon_{\infty}$ denote the composition of $\Upsilon$ with
the projection onto $\lfii(\dm)$. The restriction
$\Upsilon_{\infty,\,>1}: \widehat{\mathcal
M}_{\infty,\,>1}\to\lfii(\dm)$ is onto while for $c\in K$, the
maps $\Upsilon_{c,\,>1}$ are zero. Notice also that if $c\in
Sing(\dm)$, then
\[
\cfdm^c \cap (\oplus_{d\neq c}\cfdm^d\oplus \widehat{\mathcal
M}_{\infty,\,>1}) =\{ 0\},
\]
because tensoring both $\cfdm^c$ and $\oplus_{d\neq c}\cfdm^d$
with $\mathcal E^c$, one obtains two subspaces of $\mathcal E^c
\otimes \cfdm$ with different slopes. Also, $\widehat{\mathcal
M}_{\infty,\,>1}\cap (\oplus_c \cfdm^c)=\{0\}$. Thus the map
\[
\!\!\!(\oplus_c \Upsilon_c)\, \oplus\,
\Upsilon_{\infty,\,>1}:(\oplus_{c\,}
\cfdm^c)\,\oplus\,\widehat{\mathcal M}_{\infty,\,>1}\rightarrow \
(\bigoplus_{c\in\, Sing\dm} \lfci (\dm _{\,c}))\oplus\
\lfii(\dm_{\,\infty})
\]
is an epimorphism, and then
\begin{eqnarray*}
&\dim_{\mathcal K_{\eta^{-1}}}(\oplus_{c\,}
\cfdm^c)\,\oplus\,\widehat{\mathcal M}_{\infty,\,>1}\geqslant&\\
&\dim_{\mathcal K_{\eta^{-1}}}(\bigoplus_{c\in\, Sing\dm} \lfci
(\dm _{\,c}))\oplus\ \lfii(\dm_{\,\infty})&=\dim_{\mathcal
K_{\eta^{-1}}}\,\cfdm,
\end{eqnarray*}
the last equality by step 2. It follows that $\cfdm=(\oplus_c
\cfdm^c)\oplus \widehat{\mathcal M}_{\infty,>1}$\,,  $\Upsilon =
(\oplus_c \Upsilon_c) \oplus \Upsilon_{>1}$\, and\, $\Upsilon$ is
an isomorphism, as was to be proved. $\Box$ \vs

\newcommand{\foi}{\Phi ^{\,(\,0,\,\infty)}}
\newcommand{\at}{\mathcal A_t(1)}
\newcommand{\rr}{\mathcal R_{\eta^{-1}}(\omega)}
\newcommand{\ralg}{\mathcal A(\omega)[\eta]}

\newcommand{\eci}{\mathcal E^{\,(c, \,\infty)}}
\newcommand{\eii}{\mathcal E^{\,(\infty, \,\infty)}}
\newcommand{\eio}{\mathcal E^{\,(\,\infty,\,0)}}
\newcommand{\eoi}{\mathcal E^{\,(\,0,\,\infty)}}
\newcommand{\lfio}{\mathcal F^{\,(\,\infty,\,0)}}
\newcommand{\geveta}{\mathcal K^1_{\eta^{-1}}}
\newcommand{\gevsx}{\mathcal K^s_x}
\newcommand{\cfdmg}{\widehat{\mathcal M}_{\infty}[1]}
\newcommand{\md}{\mathcal D }
\newcommand{\sym}{\mu}

For the rest of this section we will assume that the base field
$K$ is the field $\mathbb C$ of complex numbers. In this case, one
can define subrings of the rings of formal microdifferential
operators $\lfci$ ($c\in\mathbb C \cup\{\infty\}$) adjoining
suitable convergence conditions, this is well known for
$c\in\mathbb C$, see for example \cite{Pham} or \cite{Bjo}, we
will show that an analogous definition can be given for
$c=\infty$. \vs

\noindent {\it Definition:} We denote by $\mathcal E$ the set of
formal series
\[
\sum_{i\le \,r} a_i(z)\,\eta^i \ \,,r\in \mathbb Z\ \
a_i(z)\in\mathbb C[[z]]
\]
such that
\begin{itemize}
\item[a)] There exists a $\rho_0 >0$ such that all series $a_i(z)$
are convergent in the disk $\mid z \mid <\rho_0$.
\item[b)] There exists a $0<\rho<\rho_0$ and a $\theta>0$ such that
the series
\[
\sum_{k\geqslant 0}\| a_{\,-k}(z)\|_{\rho} \frac{\theta^k}{k\,!}
\]
is convergent, where $\| a_{-k}(z)\|_{\rho}=\sup_{\mid
z\mid\leqslant\rho}\mid a_{-k}(z)\mid$.
\end{itemize}
Given $c\in \mathbb C$, we denote $\eci$ the image of $\mathcal E$
by the map
\begin{eqnarray*}
\mathcal E &\longrightarrow & \lfci \\
\sum_{i\le \,r} a_i(z)\,\eta^i &\longrightarrow &\sum_{i\le \,r}
a_i(t_c)\,\eta^i
\end{eqnarray*}
and similarly, we denote $\eii$ the image of $\mathcal E$ by the
map
\begin{eqnarray*}
\mathcal E &\longrightarrow & \lfii \\
\sum_{i\le \,r} a_i(z)\,\eta^i &\longrightarrow &\sum_{i\le \,r}
a_i(t^{-1})\,\eta^i
\end{eqnarray*}

One can prove that $\eci$ is in fact a subring of $\lfci$ and that
the division theorem (1.1) holds also for the rings $\eci$ (see
loci cit.). For $c=\infty$ these facts can be proved in a similar
way, although some modifications are needed. To illustrate them,
we prove in detail that $\eii$ is a ring using a slight variation
of the seminorms of Boutet de Monvel-Kree (\cite[Chap. 4,
\S3]{Bjo}):\vs

If $(t, \eta)$ are coordinates in $\mathbb C^2$ and $\delta>0$, we
denote $\Delta_{\,\delta}\subset \mathbb C^2$ the open subset
defined by the inequalities $\mid \eta\mid < \delta,\ \mid t \mid
>\delta^{-1}$. Clearly, if $a(z)\in \mathbb C[[z]]$ is convergent in some disk
centered at $z=0$, there exists a $\delta>0$ such that $
\sup_{\mid t\mid > \delta^{-1}}\mid a(t^{-1})\mid < \infty $,
which allows to consider the following seminorms: If
$F=\sum_{k\geqslant 0}a_{m-k}(t^{-1})\eta^{m-k}\in \eii$, put
\[
N_m(F;\delta;x):=\sum_{k,\alpha,\beta\geqslant 0}\frac{2^{-k+1}\
k\,!}{(k+\alpha)!\ (k+\beta)!}\ \|
\partial^{\alpha}_{t}\,\partial^{\beta}_{\eta}\,a_{m-k}(t^{-1})\eta^{m-k}\|_{\Delta_{\delta}}
 \ x^{2k+\alpha+\beta}
\]
where $\|
\partial^{\alpha}_{t}\,\partial^{\beta}_{\eta}\,a_{m-k}(t^{-1})\eta^{m-k}\|_{\Delta_{\delta}}:=
\sup_{(t,\eta)\in\Delta_{\delta}}\{\mid
\partial^{\alpha}_{t}\,\partial^{\beta}_{\eta}\,a_{m-k}(t^{-1})\eta^{m-k}\mid\}$
\vs

\noindent{\bf Proposition:}{\it \ If $F\in\eii$, then there exist
$\delta,\, x>0$ such that $N_m(F;\delta;x)$ is convergent.

ii) Let $a_{m-k}(z)\in \mathbb C[[t^{-1}]]$ (\,$k\geqslant 0$\,)
such that all $a_{m-k}(z)$ are convergent for $\mid z\mid<\rho$,
put $F=\sum_{k\geqslant 0}a_{m-k}(t^{-1})\eta^{m-k}$. If there
exist a $\delta,\, x>0$ such that $N_m(F;\delta;x)$ is convergent,
then $F\in\eii$.

iii) If $F,G\in\eii$, then there is a $\delta>0$ such that
$N_{ord(FG)}(F\cdot G; \delta; x)\leqslant N_{ord(F)}(F; \delta;
x)\cdot N_{ord(G)}(G; \delta; x)$, and thus $F\cdot G\in\eii$.}\vs

\noindent{\it Proof }\, (cf. \cite[Ch.4,\S 3]{Bjo}):\ i) If
$F\in\eii$, then there exist constants $A,C_1,\delta>0$ such that
\[
\|a_{m-k}(t^{-1})\eta^{m-k}\|_{\Delta_{2\delta}}\leqslant A\cdot
k\,!\cdot C_1^k \ \ \mbox{for all } k\geqslant 0
\]
From Cauchy's inequalities we get
\begin{equation}
\|\partial_t^{\alpha}\,\partial_{\eta}^{\beta}\,a_{m-k}(t^{-1})
\eta^{m-k}\|_{\Delta_{\delta}}\leqslant
\alpha\,!\cdot\beta\,!\cdot \delta^{-\alpha-\beta}\cdot 2^{\alpha}
\cdot\|a_{m-k}(t^{-1})\eta^{m-k}\|_{\Delta_{2\delta}}
\end{equation}
(notice that the factor $2^{\alpha}$ would not appear if
$\Delta_{\delta}$ were a polydisk, as it is the case for usual
microdifferential operators. However, this factor will be
harmless). Since $\alpha ! \, k!\leqslant (\alpha + k)!$\,, $\beta
! \, k!\leqslant (\beta + k)!$, putting $C_2=\max\{\sqrt{C_1},
\delta/2\}$, we get from (1) that $N_m(F;\delta;x)\leqslant A\,
\sum 2^{-k+1}(C_2 \,x)^{2k+\alpha+\beta}$, and this series is
convergent if $x<C^{-1}_2$.

Items ii) and iii) are proved exactly as for the usual
microdifferential operators, see loc. cit. $\Box$ \vs

The proof of the division theorem for the usual ring of
microdifferential operators (i.e., the ring $\eoi$), as given in
\cite{Bjo}, relies on a series of combinatorial identities and on
Cauchy's inequalities for analytic functions defined in polydisks.
As shown above, this arguments can be modified replacing these
polydisks by open sets of type $\Delta_{\delta}$, and one obtains
that the division theorem (1.1) holds also for the ring $\eii$.
\vs

\noindent {\it Definition (\cite{Ra}):}  For $s\in\mathbb R^+$, we
denote by $\gevsx$ the field $\mathbb C\{x\}_s[x^{-1}]$ of
$s$-Gevrey series on the variable $x$, this is the ring of series
$\sum_{i\geqslant 0} a_i\, x$, $a_i\in \mathbb C$ such that
$\sum_{i\geqslant 0} (a_i/ (i!)^{s})\, x^{i}$ has non-zero
convergence radius. In particular, $\mathcal K^0_x$ will denote
the field $\mathbb C\{ x\}[x^{-1}]$ of germs of meromorphic
functions. For later use, we briefly recall the behavior of Gevrey
rings and vector spaces with connection under ramification (cf.
\cite[(1.3)]{Mal1}): Let $q$ be a positive integer, $s\in\mathbb
R^+$ and set $\sigma=s/q$. The assignment $y\mapsto z^q$ defines a
morphism of fields $\pi:\mathcal K_y^s \longrightarrow\mathcal
K_z^{\sigma}$. Then:
\begin{itemize}
\item[i)] If $\mathcal V$ is a $\mathcal K_y^s$-vector space with
connection $\nabla_y$, we put $\pi^{\ast}(\mathcal V):=\mathcal
K_z^{\sigma}\otimes_{\mathcal K_y^s}\mathcal V$, endowed with the
connection $\nabla$ defined by
\[
z\nabla (\varphi \otimes v)=q\, (\varphi \otimes
y\nabla_y(v))+(z\frac{d\varphi}{dz}\otimes v).
\]
If the slopes of $\mathcal V$ are $\lambda_1,\dots,\lambda_r$,
those $\pi^{\ast}(\mathcal V)$ are
$q\,\lambda_1,\dots,q\,\lambda_r$.
\item[ii)] If $\mathcal V$ is a $\mathcal K_z^{\sigma}$-vector space with
connection $\nabla_z$, we denote $\pi_{\ast}(\mathcal V)$ the set
$\mathcal V$ regarded as a vector space over $\mathcal K_y^s$ by
restriction of scalars and endowed with the connection $\nabla:=
\frac{1}{q\,z^{q-1}}\nabla_z$. If the slopes of $\mathcal V$ are
$\lambda_1,\dots,\lambda_r$, those $\pi_{\ast}(\mathcal V)$ are
$\lambda_1/q,\dots,\lambda_r/q$ (each one repeated $q$ times).
\end{itemize}

\noindent {\it Definition:}\ If $\mathbb N$ is a $\mathbb
W_{\eta}$-module, the $s$-Gevrey germ at infinity defined by
$\mathbb N$ is the $\gevsx$-vector space $\gevsx\otimes_{\mathbb
C[\eta]}\mathbb N$, endowed with the same connection as in the
formal case (that is, $ \nabla (\alpha \otimes
n)=\partial_{\eta^{-1}}(\alpha)\otimes n - \alpha \otimes
\eta^2\partial_{\eta}n$).\vs

As in the formal case, given a holonomic $\wat$-module $\dm$ and
$c\in\mathbb C$, its microlocalization
$\eci(\dm):=\eci\otimes_{\wat}\dm$ is a $\geveta$-vector space
endowed with the connection given by left multiplication by
$\eta^2\,t$ and one defines similarly $\eii(\dm)$. We have \vs

\noindent {\bf Theorem}\ ($1$-Gevrey stationary phase): {\it Let
$\dm$ be a holonomic $\wat$-module. Then the map
\[
\Upsilon^{Gev} : \geveta\otimes_{\mathbb C[\eta]}\fdm
\longrightarrow\bigoplus_{c\in\, Sing\dm \,\cup\, \{\infty\}} \eci
(\dm)
\]
given by $\Upsilon^{Gev} (\alpha\otimes \widehat{m}) = \oplus_c \
\alpha\otimes m$ is an isomorphism of $\geveta$-vector spaces with
connection. } \vs

\noindent {\it Proof:} Again as in the formal case, the map
$\Upsilon^{Gev}$ is a morphism of $\geveta$-vector spaces with
connection, and we have to prove that it is an isomorphism. The
case of a module with punctual support is easy and left to the
reader, so we assume that $\dm=\wat /\wat \cdot P$. Then, it
follows from the division theorem for the rings $\eci$ ($c\in
Sing(\dm)\cup\{\infty\}$), that the dimension of the source and
the target of $\Upsilon^{Gev}$ are equal. So, it will be enough to
prove that the map
\[
id_{\mathbb C [[\eta^{-1}]][\eta]}\otimes
\Upsilon^{Gev}:\cfdm\longrightarrow \mathbb C
[[\eta^{-1}]][\eta]\otimes_{\geveta}(\oplus_c\,\eci(\dm))
\]
is injective. But we have a commutative diagram of $\mathbb
C[[\eta^{-1}]][\eta]$-vector spaces\vs

\begin{picture}(185,90)(-10,10)
\put(45,90){$\cfdm$}

\put(100,100){$id\otimes\Upsilon^{Gev}$}

\put(70,50){$\Upsilon$}

\put(165,50){$\mathfrak m$}

\put(115,20){$\oplus_c \lfci(\dm)$}

\put(185,90){$\mathbb C
[[\eta^{-1}]][\eta]\otimes_{\geveta}(\oplus_c\eci(\dm))$}

\put(85,95){\vector(1,0){90}}

\put(67,80){\vector(1,-1){45}}

\put(185,80){\vector(-1,-1){45}}

\end{picture}

\noindent where $\mathfrak m$ is given by $\mathfrak m (\varphi
\otimes (\oplus_c \,\xi_c))=\oplus_c \,\varphi \cdot \xi_c$. Since
$\Upsilon$ is an isomorphism (by formal stationary phase), we are
done. $\Box$ \vs

{\it Remarks:\ }i) The theorem is proved in \cite{Mal2} when
$Sing(\dm)=\{0\}$ and the singularity at infinity of $\dm$ has
slopes smaller than $+1$, by quite a different method.\vspace{3mm}

ii) It is a consequence of the above proof that the map $\mathfrak
m$ is also an isomorphism. Notice that this fails if $\dm$ is not
holonomic, for example the multiplication map $\mathbb
C[[\eta^{-1}]][\eta]\otimes_{\geveta}\eci \longrightarrow \lfci$
is clearly not onto.\vspace{3mm}

iii) Let $\dm=\wat/\wat\cdot P(t,\partial_t)$ be a holonomic
$\wat$-module such that its formal slopes at infinity of $\dm$ are
smaller than $+1$. Denote $\mathcal D^i=\mathcal K_{\eta^{-1}}^i
\langle
\partial_{\eta^{-1}}\rangle\ (i=0,1)$. For some $k\geqslant 0$ we
will have $Q= \eta^{-k} P(-\partial_{\eta},\eta)\in \mathbb
C[\eta^{-1}]\langle
\partial_{\eta^{-1}}\rangle$ (using $\partial_{\eta}=-\eta^{-2}
\partial_{\eta^{-1}}$). Consider the two
complexes of $\mathbb C$-vector spaces
\[
C_i:\ \ \ \ \ \ \ \ 0\longrightarrow \mathcal K_{\eta^{-1}}^i
\stackrel{Q}{\longrightarrow}\mathcal
K_{\eta^{-1}}^i\longrightarrow 0 \hspace{1cm} (i=0,1).
\]
There is an obvious morphism of complexes $C_0\to C_1$. Since the
formal slopes at infinity of $\dm$ are smaller than $+1$, it
follows from the results of J.P.Ramis in \cite[1.5.11, 1.5.14]{Ra}
that this morphism is a quasi-isomorphism. Since the complex $C_i$
is quasi-isomorphic to $\mathbb R \mbox{Hom}_{\mathcal
D^i}(\mathcal K^i_{\eta^{-1}}\otimes \fdm\,,\, \mathcal
K^i_{\eta^{-1}})$, by the theorem above we have a
quasi-isomorphism of solution complexes
\[
\mathbb R \mbox{Hom}_{\mathcal D^0}(\mathcal
K^0_{\eta^{-1}}\otimes \fdm\,,\, \mathcal K^0_{\eta^{-1}})\cong
\oplus_{c\in Sing \dm}\mathbb R \mbox{Hom}_{\mathcal
D^1}(\eci(\dm)\,,\, \mathcal K^1_{\eta^{-1}})
\]
That is, we can compute microlocally the ``germ at infinity" of
the solution complex of $\fdm$.\vspace{3mm}

iv) The theorem above we can be applied to study to which extent
the formal decomposition of a $\mathbb C\{x\}[x^{-1}]$-vector
space with connection given by its slopes holds at the $s$-Gevrey
level, namely one has:\vs

\noindent {\bf Theorem:}\ {\it Let $s\in \mathbb R^{+}$ and let
$\mathcal V$ be a finitely dimensional $\mathcal K^s_x$-vector
space with connection. Then there exist $\gevsx$-vector spaces
with connection \\ $\mathcal V_{<1/s}, \mathcal V_{=1/s}, \mathcal
V_{>1/s}$ of formal slopes strictly smaller than $\frac{1}{s}$
(respectively, equal to $\frac{1}{s}$, strictly greater than
$\frac{1}{s}$) and an isomorphism of $\mathcal K^s_{x}$-vector
spaces with connection
\[
\mathcal V \cong \mathcal V_{<1/s} \oplus \mathcal V_{=1/s} \oplus
\mathcal V_{>1/s}.
\]}

\noindent {\it Proof:\ } We consider first the case $s=1$. From a
theorem of Malgrange--Ramis (\cite[3.2.13]{Ra}, \cite[7.1]{Ra2},
cf. also \cite{Mal1}) on algebraization of $s$-Gevrey spaces with
connection, it follows that there is a holonomic $\wat$-module
$\mathbb N$ and an isomorphism of $\mathcal K^s_x$-vector spaces
with connection $\mathcal V \cong \mathcal K^s_{\eta^{-1}}\otimes
\mathbb N$ (the right hand side denotes the $s$-Gevrey germ at
infinity defined by $\mathbb N$, although our coordinate will be
labelled $x$ instead of $\eta^{-1}$ as done before). Let $\dm$
denote the inverse Fourier transform of $\mathbb N$. Then, putting
$\mathcal V_{<1/s}=\eoi (\dm)$, $\mathcal V_{=1/s}=\oplus_{c\in\,
Sing\dm -\{0\}} \eci (\dm)$ and $\mathcal V_{>1/s}=\eii(\dm)$ the
claimed decomposition is the one given by the $1$-Gevrey
stationary phase formula. \vs

We consider next the case $s=\frac{1}{q}$, $q$ a positive integer
(compare \cite[(2.2)]{Mal1}). Consider the map $\pi: \mathcal
K^1_y\longrightarrow\mathcal K^s_x$ given by $y\mapsto x^q$. By
the previous case, we will have
\[
\pi_{\ast}(\mathcal V)\cong\pi_{\ast}(\mathcal
V)_{<1}\oplus\pi_{\ast}(\mathcal V)_{=1}\oplus\pi_{\ast}(\mathcal
V)_{>1}.
\]
We have a surjective morphism of $\mathcal K^s_x$-vector spaces
with connection $\alpha: \pi^{\ast}(\pi_{\ast}(\mathcal V))= K^s_x
\otimes \pi_{\ast}(\mathcal V)\longrightarrow\mathcal V$ given by
$\alpha (\varphi \otimes v)=\varphi\cdot v$. Notice that all
formal slopes of $\pi^{\ast}(\pi_{\ast}(\mathcal V)_{<1})$ are
strictly smaller than $q$, those of
$\pi^{\ast}(\pi_{\ast}(\mathcal V)_{=1})$ are equal to $q$, and
those of $\pi^{\ast}(\pi_{\ast}(\mathcal V)_{>1})$ are strictly
bigger than $q$. Denote $\alpha_{<q}$ the restriction of $\alpha$
to $\pi_{\ast}(\pi^{\ast}(\mathcal V)_{<1})$ (similarly in the
cases $=q$, $>q$). Since the filtration by slopes is strict, we
have a decomposition
\[
\mathcal V \cong \alpha_{<q}(\pi_{\ast}(\pi^{\ast}(\mathcal
V)_{<1}))\oplus\alpha_{=q}(\pi_{\ast}(\pi^{\ast}(\mathcal
V)_{=1}))\oplus\alpha_{>q}(\pi_{\ast}(\pi^{\ast}(\mathcal
V)_{>1}))
\]
as desired.\vspace{3mm}

Assume now $s=\frac{p}{q}$ where $p,q\geqslant 1$ are integers. We
consider the morphism $\pi: \mathcal K^s_x\longrightarrow\mathcal
K^{1/q}_z$ given by $x\mapsto z^p$. By the previous case we have a
decomposition of $\mathcal K^{1/q}_z$- vector spaces with
connection
\[
\pi^{\ast}(\mathcal V)\cong\pi^{\ast}(\mathcal
V)_{<q}\oplus\pi^{\ast}(\mathcal V)_{=q}\oplus\pi^{\ast}(\mathcal
V)_{>q}.
\]
And an injective morphism of $\mathcal K^s_x$-vector spaces with
connection $\beta: \mathcal V
\longrightarrow\pi_{\ast}(\pi^{\ast}(\mathcal V))$ given by
$v\mapsto 1\otimes v$. Denote $\beta_{<q/p}$ the composition of
$\beta$ with the projection $\pi_{\ast}(\pi^{\ast}(\mathcal
V))\twoheadrightarrow\pi_{\ast}(\pi^{\ast}(\mathcal V)_{<q})$ and
similarly for $\beta_{=q/p}$, $\beta_{>q/p}$. Again by strictness
of the filtration by formal slopes we have a decomposition
\[
\mathcal V \cong
(\mbox{Ker}(\beta_{=q/p})\cap\mbox{Ker}(\beta_{>q/p}))\oplus
(\mbox{Ker}(\beta_{<q/p})\cap\mbox{Ker}(\beta_{>q/p}))\oplus
(\mbox{Ker}(\beta_{<q/p})\cap\mbox{Ker}(\beta_{=q/p}))
\]
which, because of the behavior of formal slopes under
$\pi_{\ast}$, is the desired one.\vs

Finally, if $s\in\mathbb R -\mathbb Q$, choose a rational number
$0<p/q<s$ such that no formal slope of $\mathcal V$ is in the
interval $[1/s, q/p]$. By the Malgrange-Ramis algebraization
theorem we can assume $\mathcal V \cong \mathcal K^s_x
\otimes\mathcal W$, where $\mathcal W$ is a $\mathcal
K_x^{p/q}$-vector space with connection. Then, the previous case
applies to $\mathcal W$ and the proof is complete. $\Box$\vs

\noindent {\bf 2. The singularity at zero of $\fdm$ and an exact
sequence of vanishing cycles} \vs

In this section we introduce one more variant of the
microlocalization functors (which should correspond to Laumon's
$(\infty,0')$ local Fourier transform). Our aim is to establish
the existence of a sequence of vanishing cycles analogous to
\cite[Theorem 10]{Kat}, \cite[proof of 3.4.2]{Lau}. In this
section we will work over the complex numbers and we will consider
only the convergent (or, more precisely $1$-Gevrey) version of the
$(\infty,0)$-microlocalization, the corresponding formal version
can be obtained just by dropping all convergence conditions. \vs

\noindent {\it Definition:\ } We denote $\eio$ the set of formal
sums
\[
P=\sum_{i\leqslant r} \ a_i(\eta)\, t^i \ \ , \ \ r\in\mathbb Z
\]
such that there exists a $\rho_0 >0$ so that all series
$a_i(\eta)$ are convergent in the disk of radius $\rho_0$ centered
at $0$, and there exists a $0<\rho<\rho_0$ and a $\theta>0$ such
that the series
\[
\sum_{k\geqslant 0}\| a_{\,-k}(\eta)\|_{\rho}
\frac{\theta^k}{k\,!}
\]
is convergent, where $\| a_{-k}(\eta)\|_{\rho}=\sup_{\mid
z\mid\leqslant\rho}\mid a_{-k}(z)\mid$. We consider in $\eio$ the
multiplication rule given by
\[
P \cdot Q = \sum_{\alpha \geqslant\, 0}\ \frac{1}{\alpha\, !}\ \
\partial_{t}^{\,\alpha}P\cdot
\partial^{\,\alpha}_{\eta} Q
\]
and the morphism of $\mathbb C$-algebras
\begin{eqnarray*}
\varphi^{(\,\infty,\,0\,)}: \wat &\longrightarrow& \eio \,,\\
t\ \ &\mapsto&\ -t\\
\partial_t\ \ &\mapsto&\ \eta
\end{eqnarray*}
which endows $\eio$ with a structure of $(\wat,\wat)$-bimodule. It
is not difficult to see that $\varphi^{(\,\infty,\,0\,)}$ is flat.
\vs

\noindent {\it Remark:\ } While the ring $\eio$ is nothing but
$\eoi$, with the r\^oles of the variables $t$ and $\eta$
interchanged, the morphism $\varphi^{(\,\infty,\,0\,)}$ is {\it
not} obtained in the same way from $\varphi^{(\,0,\,\infty\,)}$
(the morphism we considered in section 1). The morphism we get
from $\varphi^{(\,0,\,\infty\,)}$ interchanging $t$ and $\eta$
will be denoted
\begin{eqnarray*}
\mu: \mathbb W_{\eta} &\longrightarrow& \eio \,\\
\eta \ \ &\mapsto&\ \eta \\
\partial_{\eta}\ \ &\mapsto&\ t
\end{eqnarray*}

\noindent {\it Definitions:\ }i) Given a $\wat$-module $\dm$, we
put
\[
\eio(\dm):=\eio \otimes_{\wat}\dm
\]
where $\eio$ is viewed as a left $\wat$-module {\it via}
$\varphi^{(\,\infty,\,0\,)}$.

\noindent ii) Given a $\mathbb W_{\eta}$-module $\mathbb N$, we
put
\[
\sym (\mathbb N):=\eio \otimes_{\mathbb W_{\eta}}\mathbb N
\]
where now $\eio$ is viewed as a left $\mathbb W_{\eta}$-module
{\it via} $\mu$. \vs

Both $\eio(\dm)$ and $\sym(\mathbb N)$ have a structure of
$\mathbb C\{t^{-1}\}[t]$-vector spaces and of $\mathbb
C\{\eta\}\langle \partial_{\eta}\rangle$-modules, where the action
of $\partial_{\eta}$ is, by definition, given by left
multiplication by $t$. Notice that in fact $\mu(\mathbb N)$ is
nothing but $\eoi(\mathbb N)$ with the variables $t$ and $\eta$
interchanged. In section 1 we considered in the
$(0,\,\infty)$-microlocalization only the structure of $\mathbb
C\{\eta^{-1}\}[\eta]$-vector space (i.e., of $\mathbb
C\{t^{-1}\}[t]$-vector space after our interchange of variables),
while now the structure of $\mathbb C\{\eta\}\langle
\partial_{\eta}\rangle$-module will be considered as well
(and in fact it will play the major r\^ole). Notice also that
$\eio(\dm)$ depends only on the $1$-Gevrey germ at infinity
defined by $\dm$ and $\mu(\mathbb N)$ depends only on $\mathbb
N_0=\mathbb C\{\eta\}\langle
\partial_{\eta}\rangle\otimes_{\mathbb W_{\eta}}\mathbb N$.\vs

\noindent (2.1)\,{\bf Proposition}\ : {\it Let $\dm$ be a
holonomic $\wat$-module. Then the map
\begin{eqnarray*}
\Upsilon^0:\sym(\fdm) &\longrightarrow&\eio (\dm)   \\
\alpha\otimes \widehat{m}&\longrightarrow& \alpha\otimes m
\end{eqnarray*}
is an isomorphism of $\mathbb C\{\eta\}\langle
\partial_{\eta}\rangle$-modules and of
$\mathbb C\{t^{-1}\}[t]$-vector spaces.}\vs

\noindent {\it Proof:\ } It is easy to check that the map is a
morphism both of $\mathbb C\{\eta\}\langle
\partial_{\eta}\rangle$-modules and of
$\mathbb C\{t^{-1}\}[t]$-vector spaces, we have to prove that it
is an isomorphism. As for the stationary phase formulas, the
theorem reduces to the case of a Dirac $\delta$-module (then one
has $\mu(\fdm)=\eio(\dm)=0$), and the case $ \dm = \wat /\wat
\cdot P(t,\partial_t)$.\vs

In this last case, we have $\fdm = \mathbb W_{\eta}/\mathbb
W_{\eta}\cdot P(-\partial_{\eta}, \eta)$, and both the source and
the target of the map $\Upsilon^0$ are isomorphic to $\eio/\eio
\cdot P(-t,\eta)$. The map $\Upsilon^0$ composed with these
isomorphisms is the identity map, and the proposition is thus
proved. $\Box$ \vs

Let $\tau$ be a coordinate in the affine line and let $\mathbb N $
be a holonomic $\mathbb
C\{\tau\}\langle\partial_{\tau}\rangle$-module. We recall next the
formalism of solutions and microsolutions of $\mathbb N$,
following \cite{Mal2}: For $r>0$, denote by $D_r$ the disk in the
complex plane centered at $\tau=0$ and of radius $r$, by
$\widetilde{D_r^{\ast}}$ the universal covering space of
$D_r-\{0\}$, and by $\mathcal O(D_r)$ (respectively, $\mathcal O
(\widetilde{D_r^{\ast}})$) the ring of holomorphic functions on
$D_r$ (respectively, on $\widetilde{D_r^{\ast}}$). Put
$\widetilde{\mathcal C}(D_r)=\mathcal
O(\widetilde{D_r^{\ast}})/\mathcal O (D_r)$. Set $\mathcal O:=
\mbox{indlim}_{r\to 0}\,\mathcal O(D_r)$, $ \widetilde {\mathcal
O}:= \mbox{indlim}_{r\to 0}\, \widetilde{\mathcal O}(D_r)$,
$\widetilde{\mathcal C}:=\mbox{indlim}_{r\to
0}\,\widetilde{\mathcal C}(D_r)$, $\md:=\mathcal O \langle
\partial_{\tau}\rangle$. We denote by $\mathbb N\mapsto D\mathbb N$
the duality functor in the category of holonomic left
$\md$-modules, recall that if $\mathbb N = \md/\md P$, then
$D\mathbb N=\md/\md \,^{t}P$ where $^{t}P$ denotes the transposed
differential operator (see e.g. \cite[V.1]{Sa}). \vspace{3mm}

\noindent Following Malgrange, we put
\begin{enumerate}
\item[i)] $\Psi(\mathbb N):=\mbox{Hom}_{\mathcal D}(D\mathbb N,
\widetilde{\mathcal O}$) \ \ (the $\mathbb C$-vector space of
``nearby cycles" of $\mathbb N$).
\item[ii)] $\Phi(\mathbb N):=\mbox{Hom}_{\mathcal D}(D\mathbb N,
\widetilde{\mathcal C})$ \ \ (the $\mathbb C$-vector space of
microsolutions of $\mathbb N$ or ``vanishing cycles").
\end{enumerate}
Between these vector spaces there are morphisms $can: \Psi(\mathbb
N)\mapsto \Phi(\mathbb N)$ (induced by the quotient map $can:
\widetilde{\mathcal O}\rightarrow \widetilde{\mathcal C}$) and
$var: \Phi(\mathbb N)\mapsto \Psi(\mathbb N)$ (induced by the only
map $var: \widetilde{\mathcal C}\rightarrow \widetilde{\mathcal
O}$ such that $var\circ can= T-Id$, where $T$ is the monodromy on
$\widetilde{\mathcal O}$). The map $can$ is an isomorphism if
$\mathbb N \cong  \mu(\mathbb N)$, the map $var$ is an isomorphism
if $\mathbb N \cong \mathbb N[\tau^{-1}]$. The assignment $\mathbb
N\mapsto (\Psi(\mathbb N), \Phi(\mathbb N), can, var)$ is
functorial. The behavior of this spaces under localization and
microlocalization is the following:
\begin{enumerate}
\item[a)] For the localization we have $\Phi(\mathbb N[\tau^{-1}])
\cong\Psi(\mathbb N[\tau^{-1}])\cong\Psi(\mathbb N)$.
\item[b)] For the microlocalization we have $\Psi(\mu(\mathbb N))
\cong\Phi(\mu (\mathbb N))\cong\Phi(\mathbb N)$.
\end{enumerate}

For a), recall that both the kernel and cokernel of $\mathbb
N\mapsto N[\tau^{-1}]$ are a direct sum of Dirac $\delta_0$'s and
$\Psi(\delta_0)=0$ ($\delta_0=\mathbb
C\{\tau\}\langle\partial_{\tau}\rangle/(\partial_{\tau})$).
Similarly, the kernel and cokernel of $\mathbb N\mapsto\mu(\mathbb
N)$ is a direct sum of copies of the $\mathcal D$-module $\mathbb
C\{\tau\}$ (see e.g. \cite[4.11.b]{Mal1}) and $\Phi(\mathbb
C\{\tau\})=0$, from which b) follows.\vs

Given a holonomic $\mathbb W_{\tau}$-module $\dm$ we denote by
$DR(\dm)$ its De Rham complex (see e.g. \cite[I.2]{Mal2}) and by
$\mbox{Sol}(\dm)_0=\mathbb R\mbox{Hom}_{\md}(\mathcal O
\otimes_{\mathbb W_{\tau}}\dm,\, ,\mathcal O)[1]$ the stalk at
zero of its solution complex. We denote by $\mathbb
H^{\ast}_c(\mathbb A^1_{\mathbb C}, DR(\dm))$ the hypercohomology
with compact supports of the De Rham complex of $\dm$.
\vspace{3mm}

The following proposition follows essentially from results of B.
Malgrange. It shows that the De Rham cohomology with compact
supports of a holonomic $\wat$-module which has slopes at infinity
strictly smaller than $+1$ can be computed locally in terms of the
germs defined by $\dm$ at its singular points and at infinity. \vs

\noindent {\bf Proposition}\ (exact sequence of vanishing cycles):
{\it Let $\dm$ be a holonomic $\wat$-module such that all its
formal slopes at infinity are strictly smaller than $+1$. Then
there is an exact sequence of $\mathbb C$-vector spaces
\[
\hspace{-1cm} 0\rightarrow \mathbb H_c^1(\mathbb A_{\mathbb C}^1,
\mbox{DR}(\dm))\rightarrow \oplus_{c\in
Sing(\dm)}\Phi(\dm_{\,c})\rightarrow \Phi(\eio(\dm))\rightarrow
\mathbb H_c^2(\mathbb A_{\mathbb C}^1, \mbox{DR}(\dm))\rightarrow
0
\]
where $\dm_c:= \mathbb
C\{t_c\}\langle\partial_t\rangle\otimes_{\wat}\dm$.}\vspace{3mm}

\noindent{\it Proof:\ } From the exact sequence\ $
0\longrightarrow\mathcal O \longrightarrow \widetilde{\mathcal O
}\longrightarrow \widetilde{\mathcal C }\longrightarrow 0 $, \ we
get
\[
0 \rightarrow \mbox{Hom}_{\md}(D\fdm_0,\mathcal O)\rightarrow
\Psi(\fdm_{\,0})\rightarrow \Phi(\fdm_{\,0})\rightarrow
\mbox{Ext}^1_{\md}(D\fdm_0,\mathcal O)\rightarrow 0
\]
(since $\mbox{Ext}^1_{\md}(D\fdm_0,\mathcal O))=0$, see e.g.
\cite[II.3]{Mal2}). We have also quasiisomorphisms
\[
\mathbb R\Gamma_c(\mathbb A^1_{\mathbb C},
DR(\dm))\cong\mbox{Sol}\,(\widehat{D\dm})_0\,[-1]\cong\mbox{Sol}(\,D\fdm)_0\,[-1],
\]
the first one follows from \cite[VI, 2.9 and VII, 1.1]{Mal2}
(since we assume that the slopes at infinity of $\dm$ are strictly
smaller than $+1$), and the second one holds because Fourier
transform and duality commute up to the transformation given by
$t\mapsto -t, \partial_t\mapsto -\partial_t$ (see e.g.
\cite[V.2.b]{Sa}).\vs

From an element of $\oplus_{c\in Sing \dm}\Phi(\dm_{\,c})$ we get,
by the Laplace transform considered in \cite[chap. XII]{Mal2}, a
multivaluated solution of $\fdm$ defined on a half-plane in
$\mathbb C$ [loc.cit., XII, 1.2]. Under our hypothesis, the module
$\fdm$ is singular only at zero and at infinity, so this solution
can be analytically prolonged and determines univocally an element
of $\Psi(\fdm_{\,0})$. This assignment establishes an isomorphism
of complex vector spaces $ \oplus_{c\in Sing
\dm}\Phi(\dm_{\,c})\simeq\Psi(\fdm_{\,0})$. On the other hand, by
proposition (2.1) we have also an isomorphism
\[
\Phi(\fdm_{\,0})\cong\Phi(\mu(\fdm))\cong\Phi(\eio(\dm)),
\]
and the proposition follows. $\Box$ \vs

\noindent {\it Remark:} Using b) above, the long exact sequence in
the proposition can be rewritten in terms of spaces of nearby
cycles instead of spaces of microsolutions, namely one has an
exact sequence
\[
\hspace{-1cm} 0\rightarrow \mathbb H_c^1(\mathbb A_{\mathbb C}^1,
\mbox{DR}(\dm))\rightarrow \oplus_{c\in
Sing(\dm)}\Psi(\mu(\dm_{\,c}))\rightarrow
\Psi(\eio(\dm))\rightarrow \mathbb H_c^2(\mathbb A_{\mathbb C}^1,
\mbox{DR}(\dm))\rightarrow 0.
\]

\vs

\noindent {\bf 3. $p$-adic microdifferential operators of finite
order} \vs

It is proved in \cite{Mal1} (and it follows also from the
$1$-Gevrey stationary phase theorem proved in section 1) that if
$\dm$ is a holonomic $\wat$-module, singular only at zero and at
infinity, and such that the formal slopes of the singularity at
infinity are strictly smaller than $+1$, then one has an
isomorphism of $\mathcal K^1_{\eta^{-1}}$-vector spaces with
connection
\[
\mathcal K^1_{\eta^{-1}}\otimes_{\mathbb C[\eta]}\fdm\cong
\mathcal E^{(0,\infty)}(\dm)
\]
Notice that, in case the singularity at infinity of $\dm$ is
regular, these $\wat$-modules are analogous to the $\ell$-adic
canonical prolongations of Gabber and Katz, which play a major
r\^ole in G. Laumon work (\cite{Lau}, see also \cite{Kat}). So, it
might be of interest to find an analogue of Malgrange's result
when $K$ is a $p$-adic field (we denote $\mid \cdot \mid$ its
absolute value, normalized by the condition $\mid p\mid=p^{-1}$).
In such case, it seems a reasonable $p$-adic version of the ring
$\mathbb C\{\eta^{-1}\}_1$ would be the ring of power series
$\sum_{i\geqslant 0}a_i\eta^{-i}$, $a_j\in K$, such that
$\sum_{i\geqslant 0}\,i!\,a_i\,\eta^{-i}$ is convergent for $\mid
\eta^{-1} \mid <1$. If we set $\omega=\sqrt[p-1]{1/p}\in\mathbb
R$, it follows easily from the classical bounds $\omega^{k-1}<\mid
k!\mid< (k+1)\omega^{k}$ that this ring is nothing but the ring
$\mathcal A_{\eta^{-1}}(\omega)$ of power series in $\eta^{-1}$
convergent in the disk $\mid \eta^{-1} \mid <\omega$. Pursuing
this analogy, to the field $\mathcal K^1_{\eta^{-1}}$ would
correspond the ring $\mathcal A_{\eta^{-1}}(\omega)[\eta]$, which
for simplicity will be denoted $\ralg$ in the sequel. Its elements
are the Laurent series $\sum_{j\leqslant r}a_j\,\eta^{j}$ with
$a_j\in K$, $r\in\mathbb Z$, such that for all $0<\rho<\omega$,
$\limsup_{j\to -\infty}\mid a_j\mid \rho^{-j}=0$. \vs

In this section we define a ring $\foi$ of $p$-adic
microdifferential operators and a corresponding microlocalization
functor $\dm\mapsto\foi(\dm)$. We prove that if
$\dm=K[t]\langle\partial_t\rangle/K[t]\langle\partial_t\rangle\cdot
P$ is a holonomic $K[t]\langle\partial_t\rangle$-module which is
singular only at zero and infinity, the singularity at infinity
has formal slopes strictly smaller that $+1$ and the singularity
at zero is solvable at radius $1$ (\cite[8.7]{CM1}), then one has
an isomorphism of $\ralg$-modules with connection
\[
\ralg \otimes_{K[\eta]} \fdm \cong \foi(\dm).
\]
where $\fdm$ denotes now the $p$-adic Fourier transform (defined
below). This isomorphism might be regarded as a $p$-adic analogue
of the theorem of Malgrange quoted above.\vs

We assume that $K$ is a spherically complete $p$-adic field (e.g.,
a finite extension of $\mathbb Q_p$). Let $t$ be a coordinate in
the affine $K$-line. We denote by $\mathcal A_t(1)$ the ring of
power series in the variable $t$ with coefficients in $K$ which
are convergent for $\mid t\mid <1$. For all $0<\lambda<1$, the
ring $\mathcal A_t(1)$ is endowed with the norm
\[
\mid \sum_{i\geqslant 0}\,a_i\,t^i \mid_{\lambda}\, = \sup_i \{\,
\mid a_i \mid \lambda^i\,\}\in \mathbb R ^{+}.
\]

Let $r\in \mathbb Z$ be an integer, set
$\omega=\sqrt[p-1]{1/p}\in\mathbb R$. We denote by $\foi [r]$ the
set of all Laurent series $\sum_{j\leqslant r}a_{j}(t)\,\eta^j$ \,
with $a_j(t)\in \mathcal A_t(1)$, such that for all $\lambda, \rho
\in\mathbb R$ with $0<\rho<\omega\cdot\lambda<\omega$, one has
\[
\limsup_{j\to - \infty} \mid a_j(t)\mid_{\,\lambda} \,
\rho^{-j}=0.
\]
Equivalently, for all $\lambda,\rho$ in the range above there is a
$C>0$ such that for all $j\leqslant r$ one has $\mid
a_j(t)\mid_{\,\lambda}\leqslant C\cdot\rho^j$ (in fact, it is
clear that given any $0<\lambda_0<1$ is enough to check this
condition holds for $\lambda>\lambda_0$). We put $\foi = \bigcup
_{r\in\, \mathbb Z} \foi [r]$.
 \vs

\noindent {\bf Proposition:} {\it If\ $\ F=\sum f_u \,\eta^u\in
\foi$ and $\ G=\sum g_v \,\eta^v\in \foi$, then
\[
F \cdot G = \sum_{\alpha \geqslant\, 0}\ \frac{1}{\alpha\, !}\ \
\partial_{\eta}^{\,\alpha}F\cdot
\partial^{\,\alpha}_{t} G \,
\, \in \foi.
\]}

{\it Proof:\,}   Write $F\cdot G= \sum r_j(t)\,\eta^j$. We have
\begin{eqnarray*}
\mid r_j(t) \mid_{\lambda} \rho^{-j} &\leqslant
&\max_{j=u+v-\alpha} \left\{\mid \frac{1}{\alpha\, !}\
u(u-1)\dots(u-\alpha+1)\ f_u\ \frac{d^{\alpha}
g_v}{dt^{\alpha}}\mid_{\,\lambda} \rho^{-j}\right\}
\\
&\leqslant &\max_{j=u+v-\alpha} \left\{\mid
u(u-1)\dots(u-\alpha+1)\mid \ \mid f_u \mid_{\,\lambda} \ \mid g_v
\mid_{\,\lambda} \ \frac{1}{\lambda^{\alpha}}\ \rho^{-j}\right\}
\\
&\leqslant &\sup_{u} \{\mid f_u\mid_{\,\lambda}\lambda^{-u}\}\
\cdot\ \sup_{v} \{\mid g_v\mid_{\,\lambda}\lambda^{-v}\} \cdot
(\frac{\rho}{\lambda})^{\alpha}
\end{eqnarray*}
where the second inequality follows from the Cauchy inequalities.
If $j\mapsto -\infty$ then either $u\mapsto -\infty$ or $v\mapsto
-\infty$ or $\alpha\mapsto \infty$, and then we are done. $\Box$
\vs

\noindent {\it Definition:} The filtered ring $\foi$ will be
called the ring of {\it $p$-adic microdifferential operators of
finite order}. The order and the principal symbol of a
microdifferential operator are defined as in the formal case.
Notice that $\foi [0]\subset\foi$ is a filtered subring and that
one has $\ralg=K[[\eta^{-1}]][\eta]\cap \foi$.\vs

\noindent {\it Definitions:}\ If $F=\sum_{u\leqslant
m}f_{u}(t)\,\eta^u\in \mathcal A_t(1)[[\eta^{-1}]][\eta]$ and
$0<\rho<\omega\cdot\lambda<\omega$, we put
\[
\| F \|_{\lambda,\rho} = \sup_u\{\mid
f_u(t)\mid_{\lambda}\rho^{-u}\}
\]
We have $F\in\foi$ if and only if $\| F \|_{\lambda,\rho} <
\infty$ for all $\lambda, \rho$ in the range above. From the proof
of the preceding proposition follows that if
 $F,G\in\foi$, then we have $\|F\cdot
G\|_{\lambda,\rho}\leqslant \|F\|_{\lambda,\rho} \cdot\|
G\|_{\lambda,\rho}$. The subscript $\lambda,\rho$ will be omitted
if no confusion may arise. We will say that $F=\sum_{u}\,f_u\,
\eta^u \in\foi$ is {\it dominant} if there is a $\lambda_0<1$ such
that for all $\lambda_0<\lambda<1$ and
$0<\rho<\omega\cdot\lambda$, one has $\mid
f_{ord(F)}\mid_{\lambda}\rho^{-ord(F)}\,=\, \|
F\|_{\lambda,\rho}$. \vs

We want to prove a division theorem for $p$-adic microdifferential
operators of finite order. We will make implicit use of the
following lemma, its proof is elementary and left to the
reader:\vs

\noindent {\bf Lemma:}{\it\ Let $f(t)=t^m \,b(t)\in\mathcal
A_t(1)$, where $b(t)$ is invertible in $\mathcal A_t(1)$ and
$m\geqslant 0$. Then, for each $\varphi\in\mathcal A_t(1)$, there
are unique $q\in\mathcal A_t(1)$ and $r\in K[t]$ of degree smaller
or equal than $m-1$ such that $\varphi = f \cdot q +r$,\, and for
all $0<\lambda<1$,\ $\mid r\mid_{\lambda}\,\leqslant\,\mid
\varphi\mid_{\lambda}$, and\, $\mid f \mid_{\lambda}\cdot\mid
q\mid_{\lambda}\,\leqslant\,\mid \varphi\mid_{\lambda}$.} \vs

\noindent {\bf Theorem:} {\it Let $F\in\foi$ be dominant and
assume that $\sigma(F)={t}^m \,b(t)$ where $\ b(t)\in \mathcal
A_t(1)$ is invertible. Then, for all $G\in\foi$ there exist unique
$Q\in \foi$ and $R_0,\dots ,R_{m-1}\in \ralg$ such that
\[
G= Q\cdot F + {t}^{m-1}\, R_{m-1} + \dots + R_0.
\]
The remainder \ ${t}^{m-1}\, R_{m-1} + \dots + R_0$ can also be
written in a unique way in the form $S_{m-1}\,{t}^{m-1}+\dots+S_0$
with $S_i\in\ralg$.}\vs

{\it Proof:\,} It is easy to see that if $F$ is dominant so is the
product $\eta^{-ord(F)}F$ is also dominant, so we can assume $F$
is of order zero. Again, multiplying $G$ by a suitable power of
$\eta$ we may assume that $\mbox{ord}(G)=0$ as well. The existence
of a unique formal solution $Q=\sum_{j\leqslant\, 0}
q_j(t)\,\eta^j$\ , $R_i=\sum_{j\leqslant\, 0} r_{i,j}\,\eta^j\
(0\leqslant i\leqslant m-1)$ to the division problem formulated
above is well-known, the solution can be obtained as follows (cf.
\cite[Ch.4, Theorem 2.6]{Bjo}): One constructs inductively power
series $q_0,q_{-1},\dots\in\at$\ such that
\begin{eqnarray*}
G-(q_0+\dots+q_{j}\,\eta^{-j})F&=&H_{j-1}+K_{j-1}\ \ \mbox{with}\
\ H_{j-1}\in\foi[j-1] \\ \mbox{and}\ \ K_{j-1}&\in& t^{m-1}\,
\ralg + \dots + \ralg.
\end{eqnarray*}
Assume $q_0,\dots,q_{j+1}$ have already been found and put
\[
\varphi_j = g_j - \sum_{(j)}\ \frac{1}{\alpha !}\
v\,(v-1)\dots(v-\alpha+1)\ q_v \ \frac{d^{\alpha}f_u}{dt^{\alpha}}
\]
where the sum runs over those $v,u,\alpha$ with $j=v+u-\alpha$,
$\alpha\geqslant0$ and $j+1\leqslant v\leqslant 0$ (it is
understood that the product $v\,(v-1)\dots(v-\alpha+1)$ is
replaced by $1$ if $\alpha=0$). Then, the next series $q_j$ is the
quotient of the division of $\varphi_j$ by $\sigma(F)$, that is,
it is defined by the equality $\varphi_j = t^m \,b(t)\, q_j\,+r_j
$, where $q_j\in\at$ and $r_j(t)=\sum_{i=0}^{m-1}r_{i,j}t^i$ is a
polynomial of degree $m-1$ at most. The formal solution to the
division problem is given by the quotient $Q=\sum_{j\leqslant\, 0}
q_j(t)\,\eta^j$\ and the series $R_i=\sum_{j\leqslant\, 0}
r_{i,j}\,\eta^j\in\ralg$\ $(1\leqslant i\leqslant m-1)$.
\vspace{4mm}

We have to prove that this formal solution is convergent in our
sense. Fix $0<\rho <\omega\cdot \lambda<\omega$ and put
$C_{\lambda,\,\rho}=\|G\|_{\lambda,\rho}/\|F\|_{\lambda,\rho}$.
Notice that because of our hypothesis on $F$ we have
$\|F\|_{\lambda,\rho}=\mid f_0 \mid_{\lambda}$. We show next, by
descending induction on $j\leqslant 0$, that $\mid q_j
\mid\rho^{-j} \leqslant C_{\lambda,\,\rho}$. We have:
\begin{eqnarray*}
\mid \ q_j\mid_{\,\lambda}\rho^{-j} &\leqslant&\max
\left\{\,\frac{1}{\mid f_0\mid_{\lambda} }\mid\
g_j\mid_{\,\lambda}\rho^{-j},\ \right. \\ &\,& \left.
\max_{u,v,\alpha} \left\{ \frac{\mid
  v(v-1)\dots(v-\alpha+1)\mid}{\mid f_0 \mid_{\lambda}
\cdot\mid\alpha\, !\mid}\ \mid q_v\mid_{\,\lambda} \
\mid\frac{d^{\alpha}f_u}{dt^{\alpha}}\mid_{\,\lambda}
\rho^{-j}\right\}\right\}\\
&\leqslant&\max_{u,v,\alpha}
\left\{\,\frac{\|G\|}{\|F\|},\frac{1}{\lambda^{\alpha}\,\|F\|}\
\mid
q_v\mid_{\,\lambda}\ \mid  f_u \mid_{\,\lambda}\rho^{-j}\right\}\\
&=&\max_{u,v,\alpha} \left\{\,
C_{\lambda,\,\rho},\frac{1}{\lambda^{\alpha}\ \| F\|}\ (\mid
 q_v\mid_{\,\lambda}\rho^{-v})\ (\mid  f_u
\mid_{\,\lambda}\rho^{-u})\
\rho^{\alpha}\right\}\\
&\leqslant&\max\left\{ C_{\lambda,\,\rho}, \ C_{\lambda,\,\rho}\,
\cdot (\frac{\rho}{\lambda})^{\alpha}\right\}\,=\,
C_{\lambda,\,\rho}\, ,
\end{eqnarray*}
where the second inequality follows from the Cauchy inequalities.
This proves the convergence of the quotient as well as the first
inequality of norms. The convergence of the series
$R_0,\dots,R_{m-1}$ is proved similarly and it is left to the
reader. For the last statement, notice that there are unique
$S_i=\sum_{j\leqslant 0}s_{i,j}\,\eta^{j}\in K[[\eta^{-1}]][\eta]$
($0\leqslant i \leqslant m-1$), such that the remainder
${t}^{m-1}\, R_{m-1} + \dots + R_0$ can be written in the form
$S_{m-1}\,{t}^{m-1}+\dots+S_0$, in fact one has
\[
s_{i,j}=\sum_{k=0}^{m-j-1}(-1)^k \ r_{i+k, j+k}\frac{(i+k)! \cdot
(j+k)!}{k!\cdot i! \cdot j!}.
\]
From this formula follows easily that $S_i\in\ralg$ for
$i=0,\dots,m-1$. $\Box$ \vs

\noindent {\it Remark:\ } It is unreasonable to expect a division
theorem without some restriction on the divisor, for example if
$\alpha\in K$ and we take $F=1-\alpha\,\eta^{-1}\in\foi[0]$, its
formal inverse is $\sum_{i\geqslant 0} \alpha^i\eta^{-i}$, which
is not convergent in our sense for
$\mid\alpha\mid>\omega^{-1}$.\vs

We will assume that there is a $\pi \in K$ such that
$\pi^{p-1}+p=0$, which we fix from now on. Then, the morphism of
$K$-algebras defined by
\begin{eqnarray*}
\varphi^{(\,c,\,\infty\,)}: \wat &\longrightarrow& \foi \,,\\
t\ \ &\mapsto&\ t/\pi\\
\partial_t\ \ &\mapsto&\ \pi\cdot\eta
\end{eqnarray*}
endows $\foi$ with a structure of $(\wat,\wat)$-bimodule.\vs

\noindent {\it Definition:} Let $\dm$ be a $\wat$-module. We
define its $p$-adic $(0,\infty)$-microlocalization as the
$\ralg$-module
\[
\foi(\dm):= \foi \otimes_{\wat} \dm \, ,
\]
endowed with the connection given by left multiplication by
$\eta^2\cdot t$.\vs

\noindent {\it Definition} (\cite{Huy} and \cite{Meb1}): If $\dm$
is a $\wat$-module, its $p$-adic Fourier transform is defined as
$\widehat{\dm} = \mathbb W_{\eta}\otimes_{\wat } \dm$, where
$\mathbb W_{\eta}$ is regarded as a right $\wat$-module via the
$K$-algebra isomorphism given by $t\mapsto
-\,\partial_{\eta}/\pi$,\ $\partial_t \mapsto \pi\cdot \eta$. If
$m\in\dm$, we put $\widehat{m}=1\otimes m\in\widehat{\dm}$. \vs

If $\tau$ is a coordinate, we denote by $\mathcal
R_{\tau}(\theta)$ the Robba ring of power series
$\sum_{i\in\,\mathbb Z}a_i\,\tau^i$, $a_i\in K$, convergent in
some annulus $\theta-\epsilon<\,\mid\tau\mid<\theta$,
$\epsilon>0$, endowed with the derivation $\partial_{\tau}$. Let
$P(t,\partial_t)=\sum_{k=0}^d \, a_k(t)\partial_t^k\in K[t]\langle
\partial_t \rangle$, set $\mathbb M_P :=\mathbb W_t/\mathbb W_t \, P$.
We make the following assumptions on the differential operator
$P$:
\begin{enumerate}
\item $P$ is
singular only at zero and at infinity, and $\deg(a_d(t))\geqslant
\deg(a_i(t))$ for all $1\leqslant i \leqslant d$ (that is, the
formal slopes of the singularity at infinity of $P$ are smaller or
equal than $+1$).
\item The $\mathcal R_t(1)-$module with connection
$\mathcal R_t(1)\otimes_{K[t,t^{-1}]}\dm_P$ is soluble at $1$ (see
\cite[8.7]{CM1}).
\end{enumerate}
\vs

As in the formal case, if $\mathbb N$ is a $\mathbb
W_{\eta}$-module, on the $\ralg$-module $\ralg
\otimes_{K[\eta]}\mathbb N$ we will consider the connection given
by
\[
\nabla (\alpha \otimes n):=
\partial_{\eta^{-1}}(\alpha) \otimes n\, -\, \alpha \otimes
\eta^2
\partial_{\eta} n.
\]

\noindent{\bf Theorem:\ } {\it Let $P\in K[t]\langle
\partial_t \rangle$ satisfy the conditions i) and ii) above.
Then the map
\[
\Upsilon: \ralg \otimes_{K[\eta]} \widehat{\dm_P} \longrightarrow
\foi(\dm_P)
\]
given by $\Upsilon(\alpha\otimes \widehat{m})= \alpha\otimes m$ is
an isomorphism of $\ralg$-modules with connection.}\vs

\noindent {\it Proof:} It is easy to check that $\Upsilon$ is a
morphism of $\ralg$-modules with connection, we have to prove that
it is an isomorphism. We have
\[
\foi(\dm_P)\cong \frac{\foi}{\foi \cdot P(t/\pi,\pi\,\eta)}.
\]
By ii) we have (with the notations of \cite{CM1}) $Ray(\mathcal
R_t(1)\otimes_{K[t,t^{-1}]}\dm_P, 1^{-})=1$, so for $\lambda$
close to $+1$ we have $ \| a_{d-i}(t)\|_{\lambda} \leqslant
\lambda^{-i}\ \| a_d(t)\|_{\lambda}$ (\cite[Corollaire 6.4]{CM1}).
Since $P(t,\partial_t)$ is singular only at zero and infinity, we
have $a_d(t)=\alpha_d\,t^{\delta}$, with $\alpha_d\in K$, so
$\|a_d(t/\pi)\|_{\lambda}=\omega^{-\delta}\,\|a_d(t)\|_{\lambda}$
and since $\delta\geqslant \deg \, a_{d-i}(t)$ for all
$i=0,\dots,d$, we get
$\|a_{d-i}(t/\pi)\|_{\lambda}\leqslant\omega^{-\delta}\,\|a_{d-i}(t)\|_{\lambda}$,
so it follows that $\| a_{d-i}(t/\pi)\|_{\lambda} \leqslant
\lambda^{-i}\ \| a_d(t/\pi)\|_{\lambda}$. If
$0<\rho<\omega\lambda$, then
\[
\|a_d(t/\pi)\|_{\lambda} \,
\omega^{d}\rho^{-d}\geqslant\|a_{d-i}(t/\pi)\|_{\lambda} \,
\omega^{d-i}\cdot
(\omega\,\lambda)^{i}\rho^{-d}\geqslant\|a_{d-i}(t/\pi)\|_{\lambda}
\, \omega^{d-i}\cdot \rho^{-d+i}
\]
for $\lambda$ close to one, which is just the condition of
dominance for the microdifferential operator $P(t/\pi,\pi\,\eta)$.
As in the formal case, it follows now from the division theorem
that $\foi(\dm_P)$ is a free $\ralg$-module with basis $1\otimes
t^i$, $i=0,\dots,\delta-1$, and since $\Upsilon (1\otimes (-1)^i
\partial_{\eta}^i/\pi^i)=1\otimes t^i$, the morphism $\Upsilon$ is surjective. \vs

We compute next the rank over $\ralg$ of $\ralg \otimes_{\mathbb
C[\eta]} \widehat{\dm_P}$. Denote by $\alpha_i\in K$ the (possibly
zero) coefficient of $t^{\delta}$ in $a_i(t)\in K[t]$. The
coefficient of $\partial_{\eta}^{\delta}$ in the differential
operator $\widehat{P}=P(-\partial_{\eta}/\pi, \pi\,\eta)$ will be
the polynomial $q(\eta)=(-1)^{\delta}\sum_i
\alpha_{d-i}\,\pi^{d-\delta-i}\,\eta^{d-i}\in K[\eta]$. By
condition ii), $\mid\alpha_d\mid\geqslant \mid \alpha_i\mid$,
which implies that the roots of $q(\eta)$ are either zero or of
absolute value smaller than $\omega^{-1}$. It follows that
$q(\eta)$ is a unit of $\ralg$, and then the rank of $\ralg
\otimes_{\mathbb C[\eta]} \widehat{\dm_P}$ over $\ralg$ equals
$\delta$. \vs

Thus $\Upsilon$ is an epimorphism between two $\ralg$-modules with
connection of the same rank. Since the ring $\ralg$ is a
localization of $\mathcal A_{\eta^{-1}}(\omega)$, it follows from
\cite[8.1]{CM1}, that its kernel is zero and then the theorem is
proved. $\Box$\vs

Since $\ralg$ is a subring of $\rr$, we have \vs

{\bf Corollary:\ } {\it Under the same hypothesis of the theorem
above, there is an isomorphism of $\rr$-modules with connection
\[
\rr \otimes_{\mathbb K[\eta]} \widehat{\dm_P} \longrightarrow \rr
\otimes_{\ralg} \foi(\dm_P).
\]}

\noindent {\it Remark:} Given $P\in K[t]\langle
\partial_t\rangle$ verifying the conditions of the previous
theorem, one can consider as well the formal Fourier transform of
$\widehat{\dm_P}^{for}$ of $\dm_P$ (that is, the one given by
$t\mapsto -\partial_{\eta},\,
\partial_t\mapsto\eta$). One can also define a variant $\foi_1$ of
the ring $\foi$ (taking $\omega=1$), it is easy to check that one
obtains also a ring for which the division theorem holds. Then,
similarly as in the theorem above, one can show that there is an
isomorphism of $\mathcal R_{\eta^{-1}}(1)$-modules with connection
\[
\mathcal R_{\eta^{-1}}(1) \otimes_{\mathbb C[\eta]}
\widehat{\dm_P}^{for} \longrightarrow \mathcal R_{\eta^{-1}}(1)
\otimes \foi_1(\dm_P).
\]
However, unlike the $p$-adic Fourier transform, this formal
transformation does not extend to the weak completion of the Weyl
algebra considered in \cite{NM} and \cite{Huy}, and I ignore
whether it can be related to the sheaf-theoretic Fourier transform
considered in \cite{Huy}.\vs

\begin{footnotesize}

\end{footnotesize}

\

\noindent {\it Departament d'Algebra i Geometria, Universitat de
Barcelona. \\ Gran Via, 585. \\ E-08007 Barcelona, Spain. \\
e-mail adress: rgarcia@mat.ub.es}

\end{document}